\documentclass[10pt]{article}
\usepackage{amsmath}
\usepackage{amssymb}
\usepackage{amsthm}
\usepackage{mathrsfs}
\usepackage[compress,noadjust]{cite}
\numberwithin{equation}{section}

\begin{document}

\title{Endpoint estimates for commutators of sublinear operators in the Morrey type spaces}
\author{Hua Wang \footnote{E-mail address: wanghua@pku.edu.cn.}\\
\footnotesize{College of Mathematics and Econometrics, Hunan University, Changsha 410082, P. R. China}}
\date{}
\maketitle

\begin{abstract}
Let $[b,\mathcal T_\alpha]~(0\leq\alpha<n)$ be the commutators generated by $BMO(\mathbb R^n)$ functions and a class of sublinear operators satisfying certain size conditions. The aim of this paper is to study the endpoint estimates of these commutators in the weighted Morrey spaces and in the generalized Morrey spaces, under the assumptions that $[b,\mathcal T_\alpha]~(0\leq\alpha<n)$ satisfy (weighted or unweighted) endpoint inequalities on $\mathbb R^n$ and on bounded domains. Furthermore, as applications of our main results, we will obtain, in the endpoint case, the boundedness properties of many important operators in classical harmonic analysis on the weighted Morrey and the generalized Morrey spaces.
\end{abstract}

\section{Introduction and main results}

Suppose that $\mathcal T$ represents a linear or a sublinear operator, which satisfies that for any $f\in L^1(\mathbb R^n)$ with compact support and $x\notin supp\, f$,
\begin{equation}\label{sublinear}
\big|\mathcal Tf(x)\big|\leq c_0\int_{\mathbb R^n}\frac{|f(y)|}{|x-y|^n}dy,
\end{equation}
where $c_0$ is a universal constant independent of $f$ and $x\in\mathbb R^n$. The condition (\ref{sublinear}) was first introduced by Soria and Weiss in \cite{soria}. It can be proved that (\ref{sublinear}) is satisfied by many integral operators in Harmonic Analysis, such as the Hardy--Littlewood maximal operator, Calder\'on--Zygmund singular integral operators, Carleson's maximal operator, Ricci--Stein's oscillatory singular integrals and Bochner--Riesz means at the critical index and so on.

Similarly, for given $0<\alpha<n$, we assume that $\mathcal T_\alpha$ represents a linear or a sublinear operator with order $\alpha$, which satisfies that for any $f\in L^1(\mathbb R^n)$ with compact support and $x\notin supp\, f$,
\begin{equation}\label{frac sublinear}
\big|\mathcal T_\alpha f(x)\big|\leq c_1\int_{\mathbb R^n}\frac{|f(y)|}{|x-y|^{n-\alpha}}dy,
\end{equation}
where $c_1$ is also a universal constant independent of $f$ and $x\in\mathbb R^n$. It can be easily checked that (\ref{frac sublinear}) is satisfied by some important operators such as the fractional maximal operator, Riesz potential operators and fractional oscillatory singular integrals and so on.

Let $b$ be a locally integrable function on $\mathbb R^n$, suppose that the commutator operator $[b,\mathcal T]$ stands for a linear or a sublinear operator, which satisfies that for any $f\in L^1(\mathbb R^n)$ with compact support and $x\notin supp\, f$,
\begin{equation}\label{sublinear commutator}
\big|[b,\mathcal T](f)(x)\big|\leq c_2\int_{\mathbb R^n}\frac{|b(x)-b(y)|\cdot|f(y)|}{|x-y|^n}dy,
\end{equation}
where $c_2$ is an absolute constant independent of $f$ and $x\in\mathbb R^n$. Similarly, for given $0<\alpha<n$, we assume that the commutator operator $[b,\mathcal T_\alpha]$ stands for a linear or a sublinear operator, which satisfies that for any $f\in L^1(\mathbb R^n)$ with compact support and $x\notin supp\, f$,
\begin{equation}\label{frac sublinear commutator}
\big|[b,\mathcal T_\alpha](f)(x)\big|\leq c_3\int_{\mathbb R^n}\frac{|b(x)-b(y)|\cdot|f(y)|}{|x-y|^{n-\alpha}}dy,
\end{equation}
where $c_3$ is also an absolute constant independent of $f$ and $x\in\mathbb R^n$.

The classical Morrey spaces $\mathcal L^{p,\lambda}$ were originally introduced by Morrey in \cite{morrey} to study the local behavior of solutions to second order elliptic partial differential equations. Since then, these spaces play an important role in studying the regularity of solutions to partial differential equations. For the boundedness of the Hardy--Littlewood maximal operator, the fractional integral operator and the Calder\'on--Zygmund singular integral operator on these spaces, we refer the reader to \cite{adams,chiarenza,peetre}. In \cite{mizuhara}, Mizuhara introduced the generalized Morrey space $L^{p,\Theta}$ which was later extended and studied by many authors (see \cite{guliyev1,guliyev2,guliyev3,lu,nakai}). In \cite{komori}, Komori and Shirai defined the weighted Morrey space $L^{p,\kappa}(w)$ which may be
viewed as an natural generalization of weighted Lebesgue space, and then discussed the boundedness of several classical operators in Harmonic Analysis on these weighted spaces.

In \cite{lu,shi}, the authors investigated the boundedness of sublinear operators $\mathcal T_\alpha$($0\leq\alpha<n$) and their commutators with BMO functions on weighted Morrey spaces and generalized Morrey spaces. Motivated by the works in \cite{lu,shi}, in this paper, we will study the endpoint estimates of these commutators generated by $BMO(\mathbb R^n)$ functions and sublinear operators defined above in the weighted Morrey spaces $L^{1,\kappa}(w)$ for $0<\kappa<1$, and in the generalized Morrey spaces $L^{1,\Theta}$, where $\Theta$ is a growth function on $(0,+\infty)$ satisfying the doubling condition.
In order to simplify the notations, for any given $\sigma>0$, we set
\begin{equation*}
\Phi\left(\frac{|f(x)|}{\sigma}\right)=\frac{|f(x)|}{\sigma}\cdot\left(1+\log^+\frac{|f(x)|}{\sigma}\right)
\end{equation*}
when $\Phi(t)=t\cdot(1+\log^+t)$. The main results of this paper can be stated as follows.

\newtheorem{theorem}{Theorem}[section]
\newtheorem{corollary}[theorem]{Corollary}

\begin{theorem}\label{mainthm:1}
Let $0<\kappa<1$, $w\in A_1$, $b\in BMO(\mathbb R^n)$ and $[b,\mathcal T]$ satisfies the condition $(\ref{sublinear commutator})$. Suppose that for any given $\sigma>0$,
\begin{equation}\label{assump1}
w\big(\big\{x\in\mathbb R^n:\big|[b,\mathcal T](f)(x)\big|>\sigma\big\}\big)\leq C_0\int_{\mathbb R^n}\Phi\left(\frac{|f(x)|}{\sigma}\right)\cdot w(x)\,dx,
\end{equation}
where $\Phi(t)=t\cdot(1+\log^+t)$ and $C_0$ depends only on $n,w$ and $\|b\|_*$, but not on $f$ and $\sigma$. Then for the above given $\sigma>0$ and any ball $B\subset\mathbb R^n$, there exists a constant $C>0$ independent of $f$, $B$ and $\sigma$ such that
\begin{equation*}
\frac{1}{w(B)^\kappa}\cdot w\big(\big\{x\in B:\big|[b,\mathcal T](f)(x)\big|>\sigma\big\}\big)\leq C\cdot\sup_B\frac{1}{w(B)^\kappa}
\int_{B}\Phi\left(\frac{|f(x)|}{\sigma}\right)\cdot w(x)\,dx.
\end{equation*}
\end{theorem}

\begin{theorem}\label{mainthm:2}
Let $b\in BMO(\mathbb R^n)$ and $[b,\mathcal T]$ satisfies the condition $(\ref{sublinear commutator})$. Suppose that $\Theta$ satisfies $(\ref{doubling})$ with $0<D(\Theta)<2^n$, and for any given $\sigma>0$,
\begin{equation}\label{assump2}
\big|\big\{x\in\mathbb R^n:\big|[b,\mathcal T](f)(x)\big|>\sigma\big\}\big|\leq C_0\int_{\mathbb R^n}\Phi\left(\frac{|f(x)|}{\sigma}\right)dx,
\end{equation}
where $\Phi(t)=t\cdot(1+\log^+t)$ and $C_0$ depends only on $n,D(\Theta)$ and $\|b\|_*$, but not on $f$ and $\sigma$. Then for the above given $\sigma>0$ and any ball $B(x_0,r)\subset\mathbb R^n$, there exists a constant $C>0$ independent of $f$, $B(x_0,r)$ and $\sigma$ such that
\begin{equation*}
\frac{1}{\Theta(r)}\cdot\big|\big\{x\in B(x_0,r):\big|[b,\mathcal T](f)(x)\big|>\sigma\big\}\big|\leq C\cdot\sup_{r>0}\frac{1}{\Theta(r)}
\int_{B(x_0,r)}\Phi\left(\frac{|f(x)|}{\sigma}\right)dx.
\end{equation*}
\end{theorem}

\begin{theorem}\label{mainthm:3}
Let $0<\alpha<n$, $q=n/{(n-\alpha)}$, $0<\kappa<1/q$, $w^q\in A_1$, $b\in BMO(\mathbb R^n)$ and $[b,\mathcal T_\alpha]$ satisfies the condition $(\ref{frac sublinear commutator})$. Suppose that for any given $\sigma>0$ and any bounded domain $\Omega\subset\mathbb R^n$,
\begin{equation}\label{assump3}
\left[w^q\big(\big\{x\in\Omega:\big|[b,\mathcal T_\alpha](f)(x)\big|>\sigma\big\}\big)\right]^{1/q}\leq C_0\int_{\Omega}\Phi\left(\frac{|f(x)|}{\sigma}\right)\cdot w(x)\,dx,
\end{equation}
where $\Phi(t)=t\cdot(1+\log^+t)$ and $C_0$ depends only on $n,\alpha,w$ and $\|b\|_*$, but not on $f$, $\Omega$ and $\sigma$. Then for the above given $\sigma>0$ and any ball $B\subset\mathbb R^n$, there exists a constant $C>0$ independent of $f$, $B$ and $\sigma$ such that
\begin{equation*}
\begin{split}
&\left(\frac{1}{w^q(B)^{\kappa q}}\cdot w^q\big(\big\{x\in B:\big|[b,\mathcal T_{\alpha}](f)(x)\big|>\sigma\big\}\big)\right)^{1/q}\\
\leq& C\cdot\sup_B\frac{1}{w^q(B)^\kappa}
\int_{B}\Phi\left(\frac{|f(x)|}{\sigma}\right)\cdot w(x)\,dx.
\end{split}
\end{equation*}
\end{theorem}

\begin{theorem}\label{mainthm:4}
Let $0<\alpha<n$, $q=n/{(n-\alpha)}$, $b\in BMO(\mathbb R^n)$ and $[b,\mathcal T_\alpha]$ satisfies the condition $(\ref{frac sublinear commutator})$. Suppose that $\Theta$ satisfies $(\ref{doubling})$ and $0<D(\Theta)<2^{n/q}$, and for any given $\sigma>0$ and any bounded domain $\Omega\subset\mathbb R^n$,
\begin{equation}\label{assump4}
\big|\big\{x\in\Omega:\big|[b,\mathcal T_\alpha](f)(x)\big|>\sigma\big\}\big|^{1/q}\leq C_0\int_{\Omega}\Phi\left(\frac{|f(x)|}{\sigma}\right)dx,
\end{equation}
where $\Phi(t)=t\cdot(1+\log^+t)$ and $C_0$ depends only on $n,\alpha,D(\Theta)$ and $\|b\|_*$, but not on $f$, $\Omega$ and $\sigma$. Then for the above given $\sigma>0$ and any ball $B(x_0,r)\subset\mathbb R^n$, there exists a constant $C>0$ independent of $f$, $B(x_0,r)$ and $\sigma$ such that
\begin{equation*}
\begin{split}
&\left(\frac{1}{\Theta^q(r)}\cdot\big|\big\{x\in B(x_0,r):\big|[b,\mathcal T_\alpha](f)(x)\big|>\sigma\big\}\big|\right)^{1/q}\\
\leq& C\cdot\sup_{r>0}\frac{1}{\Theta(r)}\int_{B(x_0,r)}\Phi\left(\frac{|f(x)|}{\sigma}\right)dx.
\end{split}
\end{equation*}
\end{theorem}

\newtheorem{remark}[theorem]{Remark}
\begin{remark}
It should be pointed out that the conclusions of our main theorems are natural generalizations of the corresponding endpoint estimates on the weighted or unweighted Lebesgue spaces. The operators satisfying the assumptions of the above theorems include $\theta$-type Calder\'on--Zygmund operators, Marcinkiewicz integral operators, Littlewood--Paley operators, Bochner--Riesz means, fractional maximal functions and fractional integrals, which will be discussed in the last section.
\end{remark}

\section{Notations and preliminaries}

A weight $w$ will always mean a non-negative, locally integrable function on $\mathbb R^n$ which is positive on a set of positive measure, $B=B(x_0,r_B)=\{x\in\mathbb R^n:|x-x_0|<r_B\}$ denotes the open ball centered at $x_0$ and with radius $r_B>0$. Given a ball $B$ and $\lambda>0$, $\lambda B$ denotes the ball with the same center as $B$ whose radius is $\lambda$ times that of $B$. Given a Lebesgue measurable set $E$ and a weight function $w$, $|E|$ will denote the Lebesgue measure of $E$ and $w(E)=\int_E w(x)\,dx$. For $1<p<\infty$, a weight function $w$ is said to belong to the Muckenhoupt's class $A_p$, if there is a constant $C>0$ such that for every ball $B\subseteq \mathbb R^n$(see \cite{garcia,muckenhoupt}),
\begin{equation*}
\left(\frac1{|B|}\int_B w(x)\,dx\right)\left(\frac1{|B|}\int_B w(x)^{-1/{(p-1)}}\,dx\right)^{p-1}\le C.
\end{equation*}
For the case $p=1$, $w\in A_1$, if there is a constant $C>0$ such that for every ball $B\subseteq \mathbb R^n$,
\begin{equation*}
\frac1{|B|}\int_B w(x)\,dx\le C\cdot\underset{x\in B}{\mbox{ess\,inf}}\;w(x).
\end{equation*}
We also define $A_\infty=\cup_{1\leq p<\infty}A_p$. It is well known that if $w\in A_p$ with $1\leq p<\infty$, then for any ball $B$, there exists an absolute constant $C>0$ such that
\begin{equation}\label{weights}
w(2B)\le C\,w(B).
\end{equation}
In general, for $w\in A_1$ and any $\lambda>1$, there exists an absolute constant $C>0$ such that (see \cite{garcia})
\begin{equation*}
w\big(\lambda B\big)\leq C\cdot \lambda^{n}w(B).
\end{equation*}
Moreover, if $w$ is in $A_\infty$, then for all balls $B$ and all measurable subsets $E$ of $B$, there exists a number $\delta>0$ independent of $E$ and $B$ such that (see \cite{garcia})
\begin{equation}\label{compare}
\frac{w(E)}{w(B)}\leq C\left(\frac{|E|}{|B|}\right)^\delta.
\end{equation}
We say that a weight $w$ is in the reverse H\"{o}lder class $RH_s$, if there exist two constants $s>1$ and
$C>0$ such that the following reverse H\"{o}lder inequality with exponent $s>1$ holds for every
ball $B\subseteq \mathbb R^n$.
\begin{equation*}
\left(\frac{1}{|B|}\int_B w(x)^s\,dx\right)^{1/s}\le C\left(\frac{1}{|B|}\int_B w(x)\,dx\right).
\end{equation*}
Given a weight function $w$ on $\mathbb R^n$, for $1\leq p<\infty$, the weighted Lebesgue space $L^p_w(\mathbb R^n)$ is defined as the set of all functions $f$ such that
\begin{equation*}
\big\|f\big\|_{L^p_w}=\bigg(\int_{\mathbb R^n}|f(x)|^pw(x)\,dx\bigg)^{1/p}<\infty.
\end{equation*}
In particular, when $w$ equals to a constant function, we will denote $L^p_w(\mathbb R^n)$ simply by $L^p(\mathbb R^n)$.

Let $0<\kappa<1$ and $u,v$ be two weight functions on $\mathbb R^n$. Then the weighted Morrey space $L^{1,\kappa}(u,v)$ is defined by (see \cite{komori})
\begin{equation*}
L^{1,\kappa}(u,v)=\left\{f\in L^1_{loc}(u):\big\|f\big\|_{L^{1,\kappa}(u,v)}=\sup_B\frac{1}{v(B)^{\kappa}}\int_B|f(x)|u(x)\,dx<\infty\right\},
\end{equation*}
where the supremum is taken over all balls $B$ in $\mathbb R^n$. If $u=v=w$, then we set $L^{1,\kappa}(w,w)=L^{1,\kappa}(w)$.

Let $\Theta=\Theta(r)$, $r>0$, be a growth function, that is, a positive increasing function in $(0,+\infty)$ and satisfy the following doubling condition:
\begin{equation}\label{doubling}
\Theta(2r)\leq D\cdot\Theta(r), \quad \mbox{for all }\,r>0,
\end{equation}
where $D=D(\Theta)>0$ is a doubling constant independent of $r$.
The generalized Morrey space $L^{1,\Theta}(\mathbb R^n)$ is defined as the set of all locally integrable functions $f$ for which (see \cite{mizuhara})
\begin{equation*}
\sup_{r>0;B(x_0,r)}\frac{1}{\Theta(r)}\int_{B(x_0,r)}|f(x)|\,dx<\infty,
\end{equation*}
where the supremum is taken over all balls $B(x_0,r)$ in $\mathbb R^n$. From these two definitions, for given $\sigma>0$, we may rewrite the right-hand side of the inequalities in Theorems 1.1--1.4 as $\left\|\Phi\big(\frac{|f|}{\sigma}\big)\right\|_{L^{1,\kappa}(w)}$,$\left\|\Phi\big(\frac{|f|}{\sigma}\big)\right\|_{L^{1,\Theta}}$,
$\left\|\Phi\big(\frac{|f|}{\sigma}\big)\right\|_{L^{1,\kappa}(w,w^q)}$ and $\left\|\Phi\big(\frac{|f|}{\sigma}\big)\right\|_{L^{1,\Theta}}$, respectively.

We next recall some basic definitions and facts about Orlicz spaces needed for the proof of the main results. For more information on the subject, one can see \cite{rao}. A function $\Phi$ is called a Young function if it is continuous, nonnegative, convex and strictly increasing on $[0,+\infty)$ with $\Phi(0)=0$ and $\Phi(t)\to +\infty$ as $t\to +\infty$. We define the $\Phi$-average of a function $f$ over a ball $B$ by means of the following Luxemburg norm:
\begin{equation*}
\big\|f\big\|_{\Phi,B}=\inf\left\{\sigma>0:\frac{1}{|B|}\int_B\Phi\left(\frac{|f(x)|}{\sigma}\right)dx\leq1\right\}.
\end{equation*}
An equivalent norm that is often useful in calculations is as follows(see \cite{rao,perez1}):
\begin{equation}\label{equiv norm}
\big\|f\big\|_{\Phi,B}\leq \inf_{\eta>0}\left\{\eta+\frac{\eta}{|B|}\int_B\Phi\left(\frac{|f(x)|}{\eta}\right)dx\right\}\leq 2\big\|f\big\|_{\Phi,B}.
\end{equation}
Given a Young function $\Phi$, we use $\bar\Phi$ to denote the complementary Young function associated to $\Phi$. Then the following generalized H\"older's inequality holds for any given ball $B$ (see \cite{perez1,perez2}).
\begin{equation*}
\frac{1}{|B|}\int_B|f(x)\cdot g(x)|dx\leq 2\big\|f\big\|_{\Phi,B}\big\|g\big\|_{\bar\Phi,B}.
\end{equation*}
In order to deal with the weighted case, for $w\in A_\infty$, we also need to define the weighted $\Phi$-average of a function $f$ over a ball $B$ by means of the weighted Luxemburg norm:
\begin{equation*}
\big\|f\big\|_{\Phi(w),B}=\inf\left\{\sigma>0:\frac{1}{w(B)}\int_B\Phi\left(\frac{|f(x)|}{\sigma}\right)w(x)\,dx\leq1\right\}.
\end{equation*}
It can be shown that for $w\in A_\infty$(see \cite{rao,zhang}),
\begin{equation}\label{equiv norm with weight}
\big\|f\big\|_{\Phi(w),B}\approx \inf_{\eta>0}\left\{\eta+\frac{\eta}{w(B)}\int_B\Phi\left(\frac{|f(x)|}{\eta}\right)w(x)\,dx\right\},
\end{equation}
and
\begin{equation*}
\frac{1}{w(B)}\int_B|f(x)\cdot g(x)|w(x)\,dx\leq C\big\|f\big\|_{\Phi(w),B}\big\|g\big\|_{\bar\Phi(w),B}.
\end{equation*}
Here, and in what follows, $A\approx B$ means that there exist two positive constants $C_1$ and $C_2$ such that $C_1\leq\frac AB\leq C_2$.
The young function that we are going to use is $\Phi(t)=t(1+\log^+t)$ with its complementary Young function $\bar\Phi(t)\approx\exp(t)$. In the present situation, we denote
\begin{equation*}
\big\|f\big\|_{L\log L,B}=\big\|f\big\|_{\Phi,B}, \qquad \big\|g\big\|_{\exp L,B}=\big\|g\big\|_{\bar\Phi,B};
\end{equation*}
and
\begin{equation*}
\big\|f\big\|_{L\log L(w),B}=\big\|f\big\|_{\Phi(w),B}, \qquad \big\|g\big\|_{\exp L(w),B}=\big\|g\big\|_{\bar\Phi(w),B}.
\end{equation*}
By the (weighted) generalized H\"older's inequality, we have (see \cite{perez1,zhang})
\begin{equation}\label{holder}
\frac{1}{|B|}\int_B|f(x)\cdot g(x)|dx\leq 2\big\|f\big\|_{L\log L,B}\big\|g\big\|_{\exp L,B},
\end{equation}
and
\begin{equation}\label{weighted holder}
\frac{1}{w(B)}\int_B|f(x)\cdot g(x)|w(x)\,dx\leq C\big\|f\big\|_{L\log L(w),B}\big\|g\big\|_{\exp L(w),B}.
\end{equation}

Let us now recall the definition of the space of $BMO(\mathbb R^n)$ (Bounded Mean Oscillation) (see \cite{duoand,john}).
A locally integrable function $b$ is said to be in $BMO(\mathbb R^n)$, if
\begin{equation*}
\|b\|_*=\sup_{B}\frac{1}{|B|}\int_B|b(x)-b_B|\,dx<\infty,
\end{equation*}
where $b_B$ stands for the average of $b$ on $B$, i.e., $b_B=\frac{1}{|B|}\int_B b(y)\,dy$ and the supremum is taken
over all balls $B$ in $\mathbb R^n$. Modulo constants, the space $BMO(\mathbb R^n)$ is a Banach space with respect to the norm $\|\cdot\|_*$.
By the John--Nirenberg's inequality, it is not difficult to see that for any given ball $B$ (see \cite{perez1,perez2})
\begin{equation}\label{exp}
\big\|b-b_B\big\|_{\exp L,B}\leq C\|b\|_*.
\end{equation}
Furthermore, we can also prove that for any $w\in A_\infty$ and any given ball $B$ (see \cite{zhang}),
\begin{equation}\label{weighted exp}
\big\|b-b_B\big\|_{\exp L(w),B}\leq C\|b\|_*.
\end{equation}

In the sequel, the letter $C$ always denotes a positive constant which is independent of the main parameters involved, but whose value may be different from line to line. We also use $C_0,c_2,c_3$ appearing in the first section of this paper to denote certain constants. For convenience, we write $p'=p/{(p-1)}$ for given $1<p<\infty$.

\section{Proofs of Theorems \ref{mainthm:1} and \ref{mainthm:2}}

\begin{proof}[Proof of Theorem $\ref{mainthm:1}$]
Fix a ball $B=B(x_0,r_B)\subseteq\mathbb R^n$ and decompose $f=f_1+f_2$, where $f_1=f\cdot\chi_{_{2B}}$, $\chi_{_{2B}}$ denotes the characteristic function of $2B=B(x_0,2r_B)$. For any $0<\kappa<1$, $w\in A_1$ and any given $\sigma>0$, one writes
\begin{equation*}
\begin{split}
&\frac{1}{w(B)^\kappa}\cdot w\big(\big\{x\in B:\big|[b,\mathcal T](f)(x)\big|>\sigma\big\}\big)\\
\leq &\frac{1}{w(B)^\kappa}\cdot w\big(\big\{x\in B:\big|[b,\mathcal T](f_1)(x)\big|>\sigma/2\big\}\big)
+\frac{1}{w(B)^\kappa}\cdot w\big(\big\{x\in B:\big|[b,\mathcal T](f_2)(x)\big|>\sigma/2\big\}\big)\\
:=&I_1+I_2.
\end{split}
\end{equation*}
Using the condition (\ref{assump1}) and the inequality (\ref{weights}), we get
\begin{equation*}
\begin{split}
I_1&\leq C_0\cdot\frac{1}{w(B)^\kappa}\int_{\mathbb R^n}\Phi\left(\frac{|f_1(x)|}{\sigma}\right)\cdot w(x)\,dx\\
&= C_0\cdot\frac{1}{w(B)^\kappa}
\int_{2B}\Phi\left(\frac{|f(x)|}{\sigma}\right)\cdot w(x)\,dx\\
&= C_0\cdot\frac{w(2B)^\kappa}{w(B)^\kappa}\cdot\frac{1}{w(2B)^\kappa}
\int_{2B}\Phi\left(\frac{|f(x)|}{\sigma}\right)\cdot w(x)\,dx\\
&\leq C\cdot\sup_B\left\{\frac{1}{w(B)^\kappa}
\int_{B}\Phi\left(\frac{|f(x)|}{\sigma}\right)\cdot w(x)\,dx\right\}.
\end{split}
\end{equation*}
For any $x\in B$, from the definition of (\ref{sublinear commutator}), it follows that
\begin{equation*}
\begin{split}
\big|[b,\mathcal T](f_2)(x)\big|&\leq c_2\int_{\mathbb R^n}\frac{|b(x)-b(y)|\cdot|f_2(y)|}{|x-y|^n}dy\\
&\leq c_2\big|b(x)-b_{B}\big|\cdot\int_{\mathbb R^n}\frac{|f_2(y)|}{|x-y|^n}dy+c_2\int_{\mathbb R^n}\frac{|b(y)-b_{B}|\cdot|f_2(y)|}{|x-y|^n}dy\\
&:=\mu(x)+\nu(x).
\end{split}
\end{equation*}
So we have
\begin{equation*}
\begin{split}
I_2\leq&\frac{1}{w(B)^\kappa}\cdot w\big(\big\{x\in B:\mu(x)>\sigma/4\big\}\big)+\frac{1}{w(B)^\kappa}\cdot w\big(\big\{x\in B:\nu(x)>\sigma/4\big\}\big)\\
:=&I_3+I_4.
\end{split}
\end{equation*}
For the term $I_3$, for every $x\in B$, we can easily see that
\begin{equation}\label{T(f2)}
\int_{\mathbb R^n}\frac{|f_2(y)|}{|x-y|^n}dy=\int_{(2B)^c}\frac{|f(y)|}{|x-y|^{n}}dy\leq C\sum_{j=1}^\infty\frac{1}{|2^{j+1}B|}\int_{2^{j+1}B}|f(y)|\,dy.
\end{equation}
Since $w\in A_1$, then there exists a number $s>1$ such that $w\in RH_s$. Hence, by using the above pointwise estimate (\ref{T(f2)}), Chebyshev's inequality together with H\"older's inequality and John--Nirenberg's inequality (see \cite{john}), we conclude that
\begin{equation*}
\begin{split}
I_3&\leq\frac{1}{w(B)^\kappa}\cdot\frac{\,4\,}{\sigma}\int_B \mu(x)\cdot w(x)\,dx\\
&\leq C\sum_{j=1}^\infty\frac{1}{|2^{j+1}B|}\int_{2^{j+1}B}\frac{|f(y)|}{\sigma}\,dy\\
&\times\frac{1}{w(B)^\kappa}\cdot\left(\int_{B}\big|b(x)-b_{B}\big|^{s'}dx\right)^{1/{s'}}\left(\int_{B}w(x)^sdx\right)^{1/s}\\
&\leq C\sum_{j=1}^\infty\frac{1}{|2^{j+1}B|}\int_{2^{j+1}B}\frac{|f(y)|}{\sigma}\,dy\times w(B)^{1-\kappa}.
\end{split}
\end{equation*}
Furthermore, it follows directly from the $A_1$ condition and the fact $t\leq\Phi(t)=t\cdot(1+\log^+t)$ that
\begin{equation*}
\begin{split}
I_3&= C\sum_{j=1}^\infty\frac{1}{w(2^{j+1}B)}\cdot\frac{w(2^{j+1}B)}{|2^{j+1}B|}\int_{2^{j+1}B}\frac{|f(y)|}{\sigma}\,dy\times w(B)^{1-\kappa}\\
&\leq C\sum_{j=1}^\infty\frac{1}{w(2^{j+1}B)}\cdot\underset{y\in 2^{j+1}B}{\mbox{ess\,inf}}\,w(y)\int_{2^{j+1}B}\frac{|f(y)|}{\sigma}\,dy\times w(B)^{1-\kappa}\\
&\leq C\sum_{j=1}^\infty\frac{1}{w(2^{j+1}B)}\cdot\int_{2^{j+1}B}\frac{|f(y)|}{\sigma}\cdot w(y)\,dy\times w(B)^{1-\kappa}\\
&\leq C\cdot\sup_B\left\{\frac{1}{w(B)^\kappa}
\int_{B}\Phi\left(\frac{|f(y)|}{\sigma}\right)\cdot w(y)\,dy\right\}\times\sum_{j=1}^\infty\frac{w(B)^{1-\kappa}}{w(2^{j+1}B)^{1-\kappa}}.
\end{split}
\end{equation*}
Noting that $w\in A_1\subset A_\infty$, by the inequality (\ref{compare}), we get
\begin{align}\label{T<C}
\sum_{j=1}^\infty\frac{w(B)^{1-\kappa}}{w(2^{j+1}B)^{1-\kappa}}&\leq C\sum_{j=1}^\infty\left(\frac{|B|}{|2^{j+1}B|}\right)^{\delta(1-\kappa)}\notag\\
&\leq C\sum_{j=1}^\infty\left(\frac{1}{2^{(j+1)n}}\right)^{\delta(1-\kappa)}\leq C,
\end{align}
which in turn implies that
\begin{equation*}
I_3\leq C\cdot\sup_B\left\{\frac{1}{w(B)^\kappa}
\int_{B}\Phi\left(\frac{|f(y)|}{\sigma}\right)\cdot w(y)\,dy\right\}.
\end{equation*}
Similar to the proof of (\ref{T(f2)}), for all $x\in B$, we can show the following pointwise estimate as well.
\begin{equation}\label{[b,T](f2)}
\big|\nu(x)\big|\leq C\sum_{j=1}^\infty\frac{1}{|2^{j+1}B|}\int_{2^{j+1}B}\big|b(y)-b_{B}\big|\cdot\big|f(y)\big|\,dy.
\end{equation}
Applying the above pointwise estimate (\ref{[b,T](f2)}) and Chebyshev's inequality, we have
\begin{equation*}
\begin{split}
I_4&\leq\frac{1}{w(B)^\kappa}\cdot\frac{\,4\,}{\sigma}\int_B \nu(x)\cdot w(x)\,dx\\
&\leq\frac{w(B)}{w(B)^\kappa}\cdot\frac{C}{\sigma}\sum_{j=1}^\infty\frac{1}{|2^{j+1}B|}\int_{2^{j+1}B}
\big|b(y)-b_{B}\big|\cdot\big|f(y)\big|\,dy\\
&\leq\frac{w(B)}{w(B)^\kappa}\cdot\frac{C}{\sigma}\sum_{j=1}^\infty\frac{1}{|2^{j+1}B|}\int_{2^{j+1}B}
\big|b(y)-b_{2^{j+1}B}\big|\cdot\big|f(y)\big|\,dy\\
&+\frac{w(B)}{w(B)^\kappa}\cdot\frac{C}{\sigma}\sum_{j=1}^\infty\frac{1}{|2^{j+1}B|}\int_{2^{j+1}B}
\big|b_{2^{j+1}B}-b_B\big|\cdot\big|f(y)\big|\,dy\\
&:=I_5+I_6.
\end{split}
\end{equation*}
To estimate the term $I_5$, observe that for any $a,b>0$, $\Phi(a\cdot b)\leq\Phi(a)\cdot\Phi(b)$ when $\Phi(t)=t\cdot(1+\log^+t)$. We then use the generalized H\"older's inequality with weight (\ref{weighted holder}), (\ref{weighted exp}) and (\ref{equiv norm with weight}) together with (\ref{T<C}) and the $A_1$ condition to obtain
\begin{equation*}
\begin{split}
I_5&\leq \frac{C}{\sigma}\cdot w(B)^{1-\kappa}\sum_{j=1}^\infty\frac{1}{w(2^{j+1}B)}\int_{2^{j+1}B}
\big|b(y)-b_{2^{j+1}B}\big|\cdot\big|f(y)\big|w(y)\,dy\\
&\leq \frac{C}{\sigma}\cdot w(B)^{1-\kappa}\sum_{j=1}^\infty\big\|b-b_{2^{j+1}B}\big\|_{\exp L(w),2^{j+1}B}\big\|f\big\|_{L\log L(w),2^{j+1}B}\\
&\leq \frac{C\|b\|_*}{\sigma}\cdot w(B)^{1-\kappa}\sum_{j=1}^\infty
\inf_{\eta>0}\left\{\eta+\frac{\eta}{w(2^{j+1}B)}\int_{2^{j+1}B}\Phi\left(\frac{|f(z)|}{\eta}\right)w(z)\,dz\right\}\\
&\leq \frac{C\|b\|_*}{\sigma}\cdot w(B)^{1-\kappa}\sum_{j=1}^\infty
\left\{\frac{\sigma}{w(2^{j+1}B)^{1-\kappa}}+\frac{\sigma}{w(2^{j+1}B)}\int_{2^{j+1}B}\Phi\left(\frac{|f(z)|}{\sigma}\right)w(z)\,dz\right\}\\
&\leq C\|b\|_*\cdot\left[1+\sup_B\left\{\frac{1}{w(B)^\kappa}
\int_{B}\Phi\left(\frac{|f(z)|}{\sigma}\right)\cdot w(z)\,dz\right\}\right]\times\sum_{j=1}^\infty\frac{w(B)^{1-\kappa}}{w(2^{j+1}B)^{1-\kappa}}\\
&\leq C\cdot\sup_B\left\{\frac{1}{w(B)^\kappa}
\int_{B}\Phi\left(\frac{|f(z)|}{\sigma}\right)\cdot w(z)\,dz\right\}.
\end{split}
\end{equation*}
For the last term $I_6$ we proceed as follows. Since $b\in BMO(\mathbb R^n)$, then a simple calculation shows that
\begin{equation}\label{BMO1}
\big|b_{2^{j+1}B}-b_{B}\big|\leq C\cdot(j+1)\|b\|_*.
\end{equation}
Applying the inequality (\ref{BMO1}) and the facts that $w\in A_1$ and $t\leq\Phi(t)$, we get
\begin{equation*}
\begin{split}
I_6&\leq C\cdot w(B)^{1-\kappa}\sum_{j=1}^\infty(j+1)\|b\|_*\cdot\frac{1}{|2^{j+1}B|}\int_{2^{j+1}B}\frac{|f(y)|}{\sigma}\,dy\\
&\leq C\cdot w(B)^{1-\kappa}\sum_{j=1}^\infty(j+1)\|b\|_*\cdot\frac{1}{w(2^{j+1}B)}\int_{2^{j+1}B}\frac{|f(y)|}{\sigma}\cdot w(y)\,dy\\
&\leq C\cdot\sup_B\left\{\frac{1}{w(B)^\kappa}
\int_{B}\Phi\left(\frac{|f(y)|}{\sigma}\right)\cdot w(y)\,dy\right\}\times\sum_{j=1}^\infty(j+1)\cdot\frac{w(B)^{1-\kappa}}{w(2^{j+1}B)^{1-\kappa}}.
\end{split}
\end{equation*}
Since $w\in A_1\subset A_\infty$, by using the inequality (\ref{compare}) again, we have
\begin{equation*}
\begin{split}
\sum_{j=1}^\infty(j+1)\cdot\frac{w(B)^{1-\kappa}}{w(2^{j+1}B)^{1-\kappa}}&\leq C\sum_{j=1}^\infty(j+1)\cdot\left(\frac{|B|}{|2^{j+1}B|}\right)^{\delta(1-\kappa)}\\
&\leq C\sum_{j=1}^\infty(j+1)\cdot\left(\frac{1}{2^{(j+1)n}}\right)^{\delta(1-\kappa)}\leq C.
\end{split}
\end{equation*}
Therefore
\begin{equation*}
I_6\leq C\cdot\sup_B\left\{\frac{1}{w(B)^\kappa}
\int_{B}\Phi\left(\frac{|f(y)|}{\sigma}\right)\cdot w(y)\,dy\right\}.
\end{equation*}
Summarizing the above discussions, we obtain the conclusion of the theorem.
\end{proof}

\begin{proof}[Proof of Theorem $\ref{mainthm:2}$]
For any ball $B=B(x_0,r)\subseteq\mathbb R^n$ with $x_0\in\mathbb R^n$ and $r>0$, we write $f$ as $f=f_1+f_2$, where $f_1=f\cdot\chi_{_{2B}}$. Then for each fixed $\sigma>0$, we have
\begin{equation*}
\begin{split}
&\frac{1}{\Theta(r)}\cdot\big|\big\{x\in B:\big|[b,\mathcal T](f)(x)\big|>\sigma\big\}\big|\\
\leq &\frac{1}{\Theta(r)}\cdot \big|\big\{x\in B:\big|[b,\mathcal T](f_1)(x)\big|>\sigma/2\big\}\big|
+\frac{1}{\Theta(r)}\cdot \big|\big\{x\in B:\big|[b,\mathcal T](f_2)(x)\big|>\sigma/2\big\}\big|\\
:=&J_1+J_2.
\end{split}
\end{equation*}
We consider the term $J_1$ first. The condition (\ref{assump2}) and the inequality (\ref{doubling}) imply that
\begin{equation*}
\begin{split}
J_1&\leq C_0\cdot\frac{1}{\Theta(r)}\int_{\mathbb R^n}\Phi\left(\frac{|f_1(x)|}{\sigma}\right)dx\\
&= C_0\cdot\frac{1}{\Theta(r)}
\int_{2B}\Phi\left(\frac{|f(x)|}{\sigma}\right)dx\\
\end{split}
\end{equation*}
\begin{equation*}
\begin{split}
&= C_0\cdot\frac{\Theta(2r)}{\Theta(r)}\cdot\frac{1}{\Theta(2r)}
\int_{B(x_0,2r)}\Phi\left(\frac{|f(x)|}{\sigma}\right)dx\\
&\leq C\cdot\sup_{r>0;B(x_0,r)}\left\{\frac{1}{\Theta(r)}
\int_{B(x_0,r)}\Phi\left(\frac{|f(x)|}{\sigma}\right)dx\right\}.
\end{split}
\end{equation*}
We now turn our attention to the estimate of $J_2$. Recalling that the following estimate holds for any $x\in B$,
\begin{equation*}
\big|[b,\mathcal T](f_2)(x)\big|\leq\mu(x)+\nu(x),
\end{equation*}
where
\begin{equation*}
\mu(x)=c_2\big|b(x)-b_{B}\big|\cdot\int_{\mathbb R^n}\frac{|f_2(y)|}{|x-y|^n}dy,
\end{equation*}
and
\begin{equation*}
\nu(x)=c_2\int_{\mathbb R^n}\frac{|b(y)-b_{B}|\cdot|f_2(y)|}{|x-y|^n}dy.
\end{equation*}
Thus, we have
\begin{equation*}
\begin{split}
J_2\leq&\frac{1}{\Theta(r)}\cdot \big|\big\{x\in B:\mu(x)>\sigma/4\big\}\big|+\frac{1}{\Theta(r)}\cdot \big|\big\{x\in B:\nu(x)>\sigma/4\big\}\big|\\
:=&J_3+J_4.
\end{split}
\end{equation*}
By using the previous pointwise estimate (\ref{T(f2)}), Chebyshev's inequality and the definition of BMO, we can deduce that
\begin{equation*}
\begin{split}
J_3&\leq\frac{1}{\Theta(r)}\cdot\frac{\,4\,}{\sigma}\int_B \mu(x)\,dx\\
&\leq C\sum_{j=1}^\infty\frac{1}{|2^{j+1}B|}\int_{2^{j+1}B}\frac{|f(y)|}{\sigma}\,dy
\times\left\{\frac{|B|}{\Theta(r)}\cdot\frac{1}{|B|}\int_B\big|b(x)-b_{B}\big|\,dx\right\}\\
&\leq C\|b\|_*\sum_{j=1}^\infty\frac{|B|}{|2^{j+1}B|}\cdot\frac{\Theta(2^{j+1}r)}{\Theta(r)}\cdot
\frac{1}{\Theta(2^{j+1}r)}\int_{B(x_0,2^{j+1}r)}\frac{|f(y)|}{\sigma}\,dy\\
&\leq C\cdot\sup_{r>0;B(x_0,r)}\left\{\frac{1}{\Theta(r)}
\int_{B(x_0,r)}\Phi\left(\frac{|f(y)|}{\sigma}\right)dy\right\}\times\sum_{j=1}^\infty\frac{|B|}{|2^{j+1}B|}\cdot\frac{\Theta(2^{j+1}r)}{\Theta(r)}.
\end{split}
\end{equation*}
Note that $0<D(\Theta)<2^n$, then by using the doubling condition (\ref{doubling}) of $\Theta$, we can see that
\begin{equation}\label{theta<C1}
\sum_{j=1}^\infty\frac{|B|}{|2^{j+1}B|}\cdot\frac{\Theta(2^{j+1}r)}{\Theta(r)}\leq\sum_{j=1}^\infty\left(\frac{D(\Theta)}{2^n}\right)^{j+1}\leq C,
\end{equation}
which in turn gives that
\begin{equation*}
J_3\leq C\cdot\sup_{r>0;B(x_0,r)}\left\{\frac{1}{\Theta(r)}
\int_{B(x_0,r)}\Phi\left(\frac{|f(y)|}{\sigma}\right)dy\right\}.
\end{equation*}
Applying the previous pointwise estimate (\ref{[b,T](f2)}) and Chebyshev's inequality, we have
\begin{equation*}
\begin{split}
J_4&\leq\frac{1}{\Theta(r)}\cdot\frac{\,4\,}{\sigma}\int_B \nu(x)\,dx\\
&\leq\frac{|B|}{\Theta(r)}\cdot\frac{C}{\sigma}\sum_{j=1}^\infty\frac{1}{|2^{j+1}B|}\int_{2^{j+1}B}
\big|b(y)-b_{B}\big|\cdot\big|f(y)\big|\,dy\\
&\leq\frac{|B|}{\Theta(r)}\cdot\frac{C}{\sigma}\sum_{j=1}^\infty\frac{1}{|2^{j+1}B|}\int_{2^{j+1}B}
\big|b(y)-b_{2^{j+1}B}\big|\cdot\big|f(y)\big|\,dy\\
&+\frac{|B|}{\Theta(r)}\cdot\frac{C}{\sigma}\sum_{j=1}^\infty\frac{1}{|2^{j+1}B|}\int_{2^{j+1}B}
\big|b_{2^{j+1}B}-b_B\big|\cdot\big|f(y)\big|\,dy\\
&:=J_5+J_6.
\end{split}
\end{equation*}
For the term $J_5$, notice that the inequality $\Phi(a\cdot b)\leq\Phi(a)\cdot\Phi(b)$ holds for any $a,b>0$, when $\Phi(t)=t\cdot(1+\log^+t)$. We then use the generalized H\"older's inequality (\ref{holder}), (\ref{exp}) and (\ref{equiv norm}) together with (\ref{theta<C1}) to obtain
\begin{equation*}
\begin{split}
J_5&\leq \frac{|B|}{\Theta(r)}\cdot\frac{C}{\sigma}\sum_{j=1}^\infty\big\|b-b_{2^{j+1}B}\big\|_{\exp L,2^{j+1}B}\big\|f\big\|_{L\log L,2^{j+1}B}\\
&\leq \frac{C\|b\|_*}{\sigma}\cdot\frac{|B|}{\Theta(r)}\sum_{j=1}^\infty
\inf_{\eta>0}\left\{\eta+\frac{\eta}{|2^{j+1}B|}\int_{2^{j+1}B}\Phi\left(\frac{|f(z)|}{\eta}\right)dz\right\}\\
&\leq \frac{C\|b\|_*}{\sigma}\cdot \frac{|B|}{\Theta(r)}\sum_{j=1}^\infty
\left\{\frac{\sigma\cdot\Theta(2^{j+1}r)}{|2^{j+1}B|}+\frac{\sigma}{|2^{j+1}B|}\int_{2^{j+1}B}\Phi\left(\frac{|f(z)|}{\sigma}\right)dz\right\}\\
&\leq C\|b\|_*\cdot\left[1+\sup_{r>0;B(x_0,r)}\left\{\frac{1}{\Theta(r)}
\int_{B(x_0,r)}\Phi\left(\frac{|f(z)|}{\sigma}\right)dz\right\}\right]\\
&\times\sum_{j=1}^\infty\frac{|B|}{|2^{j+1}B|}\cdot\frac{\Theta(2^{j+1}r)}{\Theta(r)}\\
&\leq C\cdot\sup_{r>0;B(x_0,r)}\left\{\frac{1}{\Theta(r)}
\int_{B(x_0,r)}\Phi\left(\frac{|f(z)|}{\sigma}\right)dz\right\}.
\end{split}
\end{equation*}
For the last term $J_6$, an application of the inequality (\ref{BMO1}) leads to that
\begin{equation*}
\begin{split}
J_6&\leq C\cdot\frac{|B|}{\Theta(r)}\sum_{j=1}^\infty(j+1)\|b\|_*\cdot\frac{1}{|2^{j+1}B|}\int_{2^{j+1}B}\frac{|f(y)|}{\sigma}\,dy\\
&= C\cdot\frac{|B|}{\Theta(r)}\sum_{j=1}^\infty(j+1)\|b\|_*\cdot\frac{\Theta(2^{j+1}r)}{|2^{j+1}B|}\cdot
\frac{1}{\Theta(2^{j+1}r)}\int_{B(x_0,2^{j+1}r)}\frac{|f(y)|}{\sigma}\,dy\\
\end{split}
\end{equation*}
\begin{equation*}
\begin{split}
&\leq C\cdot\sup_{r>0;B(x_0,r)}\left\{\frac{1}{\Theta(r)}
\int_{B(x_0,r)}\Phi\left(\frac{|f(y)|}{\sigma}\right)dy\right\}\\
&\times\sum_{j=1}^\infty(j+1)\cdot\frac{|B|}{|2^{j+1}B|}\cdot\frac{\Theta(2^{j+1}r)}{\Theta(r)}.
\end{split}
\end{equation*}
Moreover, by using the doubling condition (\ref{doubling}) of $\Theta$ again and the fact that $0<D(\Theta)<2^n$, we find that
\begin{equation}\label{theta<C2}
\sum_{j=1}^\infty(j+1)\cdot\frac{|B|}{|2^{j+1}B|}\cdot\frac{\Theta(2^{j+1}r)}{\Theta(r)}\leq C\sum_{j=1}^\infty(j+1)\cdot\left(\frac{D(\Theta)}{2^n}\right)^{j+1}\leq C.
\end{equation}
Substituting the above inequality (\ref{theta<C2}) into the term $J_6$, we thus obtain
\begin{equation*}
J_6\leq C\cdot\sup_{r>0;B(x_0,r)}\left\{\frac{1}{\Theta(r)}
\int_{B(x_0,r)}\Phi\left(\frac{|f(y)|}{\sigma}\right)dy\right\}.
\end{equation*}
Summing up all the above estimates, we therefore conclude the proof of the main theorem.
\end{proof}

\section{Proofs of Theorems \ref{mainthm:3} and \ref{mainthm:4}}

\begin{proof}[Proof of Theorem $\ref{mainthm:3}$]
Fix a ball $B=B(x_0,r_B)\subseteq\mathbb R^n$ and $x\in B$, we split $f$ as usual by $f=f\cdot\chi_{_{2B}}+f\cdot\chi_{_{(2B)^c}}:=f_1+f_2$. For any $0<\kappa<1/q$, $w^q\in A_1$ with $q=n/{(n-\alpha)}>1$ and any given $\sigma>0$, we then write
\begin{equation*}
\begin{split}
&\left(\frac{1}{w^q(B)^{\kappa q}}\cdot w^q\big(\big\{x\in B:\big|[b,\mathcal T_{\alpha}](f)(x)\big|>\sigma\big\}\big)\right)^{1/q}\\
\leq &\left(\frac{1}{w^q(B)^{\kappa q}}\cdot w^q\big(\big\{x\in B:\big|[b,\mathcal T_{\alpha}](f_1)(x)\big|>\sigma/2\big\}\big)\right)^{1/q}\\
&+\left(\frac{1}{w^q(B)^{\kappa q}}\cdot w^q\big(\big\{x\in B:\big|[b,\mathcal T_{\alpha}](f_2)(x)\big|>\sigma/2\big\}\big)\right)^{1/q}\\
:=&I'_1+I'_2.
\end{split}
\end{equation*}
By using the assumption (\ref{assump3}) and the inequality (\ref{weights}), we get
\begin{equation*}
\begin{split}
I'_1&\leq C_0\cdot\frac{1}{w^q(B)^\kappa}\int_{B}\Phi\left(\frac{|f_1(x)|}{\sigma}\right)\cdot w(x)\,dx\\
&= C_0\cdot\frac{1}{w^q(B)^\kappa}
\int_{2B}\Phi\left(\frac{|f(x)|}{\sigma}\right)\cdot w(x)\,dx\\
&= C_0\cdot\frac{w^q(2B)^\kappa}{w^q(B)^\kappa}\cdot\frac{1}{w^q(2B)^\kappa}
\int_{2B}\Phi\left(\frac{|f(x)|}{\sigma}\right)\cdot w(x)\,dx\\
&\leq C\cdot\sup_B\left\{\frac{1}{w^q(B)^\kappa}
\int_{B}\Phi\left(\frac{|f(x)|}{\sigma}\right)\cdot w(x)\,dx\right\}.
\end{split}
\end{equation*}
For any $x\in B$, from the definition of (\ref{frac sublinear commutator}), it follows that
\begin{equation*}
\begin{split}
\big|[b,\mathcal T_\alpha](f_2)(x)\big|&\leq c_3\int_{\mathbb R^n}\frac{|b(x)-b(y)|\cdot|f_2(y)|}{|x-y|^{n-\alpha}}dy\\
&\leq c_3\big|b(x)-b_{B}\big|\cdot\int_{\mathbb R^n}\frac{|f_2(y)|}{|x-y|^{n-\alpha}}dy+c_3\int_{\mathbb R^n}\frac{|b(y)-b_{B}|\cdot|f_2(y)|}{|x-y|^{n-\alpha}}dy\\
&:=\widetilde\mu(x)+\widetilde\nu(x).
\end{split}
\end{equation*}
So we can rewrite the term $I'_2$ as follows:
\begin{equation*}
\begin{split}
I'_2\leq&\left(\frac{1}{w^q(B)^{\kappa q}}\cdot w^q\big(\big\{x\in B:\widetilde\mu(x)>\sigma/4\big\}\big)\right)^{1/q}\\
&+\left(\frac{1}{w^q(B)^{\kappa q}}\cdot w^q\big(\big\{x\in B:\widetilde\nu(x)>\sigma/4\big\}\big)\right)^{1/q}\\
:=&I'_3+I'_4.
\end{split}
\end{equation*}
For the term $I'_3$, for given $0<\alpha<n$ and every $x\in B$, we can easily check that
\begin{equation}\label{Ta(f2)}
\begin{split}
\int_{\mathbb R^n}\frac{|f_2(y)|}{|x-y|^{n-\alpha}}dy&=\int_{(2B)^c}\frac{|f(y)|}{|x-y|^{n-\alpha}}dy\\
&\leq C\sum_{j=1}^\infty\frac{1}{|2^{j+1}B|^{1-\alpha/n}}\int_{2^{j+1}B}|f(y)|\,dy.
\end{split}
\end{equation}
Since $w^q$ is in $A_1$, we know that there exists a number $r>1$ such that $w^q\in RH_r$. Hence, by using the above pointwise estimate (\ref{Ta(f2)}), Chebyshev's inequality together with H\"older's inequality and John--Nirenberg's inequality (see \cite{john}), we deduce that
\begin{equation*}
\begin{split}
I'_3&\leq\frac{1}{w^q(B)^\kappa}\cdot\frac{\,4\,}{\sigma}\left(\int_B|\widetilde\mu(x)|^qw^q(x)\,dx\right)^{1/q}\\
&\leq C\sum_{j=1}^\infty\frac{1}{|2^{j+1}B|^{1-\alpha/n}}\int_{2^{j+1}B}\frac{|f(y)|}{\sigma}\,dy
\times\frac{1}{w^q(B)^\kappa}\cdot\left(\int_B\big|b(x)-b_{B}\big|^qw^q(x)\,dx\right)^{1/q}\\
&\leq C\sum_{j=1}^\infty\frac{1}{|2^{j+1}B|^{1-\alpha/n}}\int_{2^{j+1}B}\frac{|f(y)|}{\sigma}\,dy\\
&\times\frac{1}{w^q(B)^\kappa}\cdot\left(\int_{B}\big|b(x)-b_{B}\big|^{qr'}dx\right)^{1/{(qr')}}\left(\int_{B}\big[w^q(x)\big]^rdx\right)^{1/{(qr)}}\\
&\leq C\sum_{j=1}^\infty\frac{1}{|2^{j+1}B|^{1-\alpha/n}}\int_{2^{j+1}B}\frac{|f(y)|}{\sigma}\,dy\times w^q(B)^{1/q-\kappa}.\\
\end{split}
\end{equation*}
Moreover, by applying H\"older's inequality and then the reverse H\"older inequality in succession, we can show that $w\in A_1\cap RH_q$ if and only if $w^q\in A_1$(see \cite{johnson}). Thus, we are able to verify that for any $j\in\mathbb Z_+$,
\begin{equation*}
w^q(2^{j+1}B)^{1/q}=\left(\int_{2^{j+1}B}w^q(x)\,dx\right)^{1/q}\leq C\cdot|2^{j+1}B|^{1/q-1}\cdot w(2^{j+1}B),
\end{equation*}
which is equivalent to
\begin{equation}\label{wq(B)/B}
\frac{w^q(2^{j+1}B)^{1/q}}{|2^{j+1}B|^{1/q}}\leq C\cdot\frac{w(2^{j+1}B)}{|2^{j+1}B|}.
\end{equation}
Therefore, by using the inequality (\ref{wq(B)/B}) together with the facts that $1/q=1-\alpha/n$, $w\in A_1$ and $t\leq\Phi(t)$, we obtain
\begin{equation*}
\begin{split}
I'_3&= C\sum_{j=1}^\infty\frac{w^q(2^{j+1}B)^{1/q}}{|2^{j+1}B|^{1-\alpha/n}}\cdot\frac{1}{w^q(2^{j+1}B)^{\kappa}}
\int_{2^{j+1}B}\frac{|f(y)|}{\sigma}\,dy\times\frac{w^q(B)^{1/q-\kappa}}{w^q(2^{j+1}B)^{1/q-\kappa}}\\
&\leq C\sum_{j=1}^\infty\frac{w(2^{j+1}B)}{|2^{j+1}B|}\cdot\frac{1}{w^q(2^{j+1}B)^{\kappa}}
\int_{2^{j+1}B}\frac{|f(y)|}{\sigma}\,dy\times\frac{w^q(B)^{1/q-\kappa}}{w^q(2^{j+1}B)^{1/q-\kappa}}\\
&\leq C\sum_{j=1}^\infty\frac{1}{w^q(2^{j+1}B)^{\kappa}}
\int_{2^{j+1}B}\frac{|f(y)|}{\sigma}\cdot w(y)\,dy\times\frac{w^q(B)^{1/q-\kappa}}{w^q(2^{j+1}B)^{1/q-\kappa}}\\
&\leq C\cdot\sup_B\left\{\frac{1}{w^q(B)^\kappa}
\int_{B}\Phi\left(\frac{|f(y)|}{\sigma}\right)\cdot w(y)\,dy\right\}\times\sum_{j=1}^\infty\frac{w^q(B)^{1/q-\kappa}}{w^q(2^{j+1}B)^{1/q-\kappa}}.
\end{split}
\end{equation*}
Applying the inequality (\ref{compare}) and the property $w^q\in A_1\subset A_\infty$, we can get
\begin{align}\label{wq(B)<C}
\sum_{j=1}^\infty\frac{w^q(B)^{1/q-\kappa}}{w^q(2^{j+1}B)^{1/q-\kappa}}
&\leq C\sum_{j=1}^\infty\left(\frac{|B|}{|2^{j+1}B|}\right)^{\delta^*(1/q-\kappa)}\notag\\
&\leq C\sum_{j=1}^\infty\left(\frac{1}{2^{(j+1)n}}\right)^{\delta^*(1/q-\kappa)}\leq C,
\end{align}
where in the last inequality we have used the facts that $\delta^*>0$ and $0<\kappa<1/q$. Substituting the above inequality (\ref{wq(B)<C}) into the term $I'_3$, we thus obtain
\begin{equation*}
I'_3\leq C\cdot\sup_B\left\{\frac{1}{w^q(B)^\kappa}
\int_{B}\Phi\left(\frac{|f(y)|}{\sigma}\right)\cdot w(y)\,dy\right\}.
\end{equation*}
For the term $I'_4$, similar to the proof of (\ref{Ta(f2)}), for all $0<\alpha<n$ and all $x\in B$, we can show the following pointwise estimate as well.
\begin{equation}\label{[b,Ta](f2)}
\big|\widetilde\nu(x)\big|\leq C\sum_{j=1}^\infty\frac{1}{|2^{j+1}B|^{1-\alpha/n}}\int_{2^{j+1}B}\big|b(y)-b_{B}\big|\cdot\big|f(y)\big|\,dy.
\end{equation}
Following the same arguments as in the proof of Theorem \ref{mainthm:1} and using the pointwise inequality (\ref{[b,Ta](f2)}) and Chebyshev's inequality, we have
\begin{equation*}
\begin{split}
I'_4&\leq\frac{1}{w^q(B)^\kappa}\cdot\frac{\,4\,}{\sigma}\left(\int_B|\widetilde\nu(x)|^qw^q(x)\,dx\right)^{1/q}\\
&\leq\frac{w^q(B)^{1/q}}{w^q(B)^\kappa}\cdot\frac{C}{\sigma}\sum_{j=1}^\infty\frac{1}{|2^{j+1}B|^{1-\alpha/n}}\int_{2^{j+1}B}
\big|b(y)-b_{B}\big|\cdot\big|f(y)\big|\,dy\\
&\leq\frac{w^q(B)^{1/q}}{w^q(B)^\kappa}\cdot\frac{C}{\sigma}\sum_{j=1}^\infty\frac{1}{|2^{j+1}B|^{1-\alpha/n}}\int_{2^{j+1}B}
\big|b(y)-b_{2^{j+1}B}\big|\cdot\big|f(y)\big|\,dy\\
&+\frac{w^q(B)^{1/q}}{w^q(B)^\kappa}\cdot\frac{C}{\sigma}\sum_{j=1}^\infty\frac{1}{|2^{j+1}B|^{1-\alpha/n}}\int_{2^{j+1}B}
\big|b_{2^{j+1}B}-b_B\big|\cdot\big|f(y)\big|\,dy\\
&:=I'_5+I'_6.
\end{split}
\end{equation*}
To deal with the term $I'_5$, it then follows from the inequality (\ref{wq(B)/B}) and the facts $1/q=1-\alpha/n$ and $w\in A_1$ that
\begin{equation*}
\begin{split}
I'_5&=\frac{C}{\sigma}\sum_{j=1}^\infty\frac{w^q(B)^{1/q-\kappa}}{w^q(2^{j+1}B)^{1/q-\kappa}}\cdot\frac{1}{w^q(2^{j+1}B)^\kappa}\\
&\times\frac{w^q(2^{j+1}B)^{1/q}}{|2^{j+1}B|^{1-\alpha/n}}\cdot\int_{2^{j+1}B}
\big|b(y)-b_{2^{j+1}B}\big|\cdot\big|f(y)\big|\,dy\\
&\leq \frac{C}{\sigma}\sum_{j=1}^\infty\frac{w^q(B)^{1/q-\kappa}}{w^q(2^{j+1}B)^{1/q-\kappa}}\cdot\frac{1}{w^q(2^{j+1}B)^\kappa}\\
&\times\frac{w(2^{j+1}B)}{|2^{j+1}B|}\cdot\int_{2^{j+1}B}
\big|b(y)-b_{2^{j+1}B}\big|\cdot\big|f(y)\big|\,dy\\
&\leq \frac{C}{\sigma}\sum_{j=1}^\infty\frac{w^q(B)^{1/q-\kappa}}{w^q(2^{j+1}B)^{1/q-\kappa}}\cdot\frac{1}{w^q(2^{j+1}B)^\kappa}\times\int_{2^{j+1}B}
\big|b(y)-b_{2^{j+1}B}\big|\cdot\big|f(y)\big|w(y)\,dy.\\
\end{split}
\end{equation*}
Furthermore, by using the generalized H\"older's inequality with weight (\ref{weighted holder}), (\ref{weighted exp}), (\ref{equiv norm with weight}) together with (\ref{wq(B)<C}) and the fact that $\Phi(t)$ is submultiplicative ($\Phi(a\cdot b)\leq\Phi(a)\cdot\Phi(b)$ for any $a,b>0$), we can conclude that
\begin{equation*}
\begin{split}
I'_5&\leq \frac{C}{\sigma}\sum_{j=1}^\infty\frac{w^q(B)^{1/q-\kappa}}{w^q(2^{j+1}B)^{1/q-\kappa}}\cdot
\frac{w(2^{j+1}B)}{w^q(2^{j+1}B)^\kappa}\big\|b-b_{2^{j+1}B}\big\|_{\exp L(w),2^{j+1}B}\big\|f\big\|_{L\log L(w),2^{j+1}B}\\
&\leq \frac{C\|b\|_*}{\sigma}\cdot\sum_{j=1}^\infty\frac{w^q(B)^{1/q-\kappa}}{w^q(2^{j+1}B)^{1/q-\kappa}}\cdot\frac{w(2^{j+1}B)}{w^q(2^{j+1}B)^\kappa}\\
&\times\inf_{\eta>0}\left\{\eta+\frac{\eta}{w(2^{j+1}B)}\int_{2^{j+1}B}\Phi\left(\frac{|f(z)|}{\eta}\right)w(z)\,dz\right\}\\
\end{split}
\end{equation*}
\begin{equation*}
\begin{split}
&\leq \frac{C\|b\|_*}{\sigma}\cdot\sum_{j=1}^\infty\frac{w^q(B)^{1/q-\kappa}}{w^q(2^{j+1}B)^{1/q-\kappa}}\cdot\frac{w(2^{j+1}B)}{w^q(2^{j+1}B)^\kappa}\\
&\times\left\{\frac{\sigma\cdot w^q(2^{j+1}B)^\kappa}{w(2^{j+1}B)}+\frac{\sigma}{w(2^{j+1}B)}\int_{2^{j+1}B}\Phi\left(\frac{|f(z)|}{\sigma}\right)w(z)\,dz\right\}\\
&\leq C\|b\|_*\cdot\left[1+\sup_B\left\{\frac{1}{w^q(B)^\kappa}
\int_{B}\Phi\left(\frac{|f(z)|}{\sigma}\right)\cdot w(z)\,dz\right\}\right]\times\sum_{j=1}^\infty\frac{w^q(B)^{1/q-\kappa}}{w^q(2^{j+1}B)^{1/q-\kappa}}\\
&\leq C\cdot\sup_B\left\{\frac{1}{w^q(B)^\kappa}
\int_{B}\Phi\left(\frac{|f(z)|}{\sigma}\right)\cdot w(z)\,dz\right\}.
\end{split}
\end{equation*}
For the last term $I'_6$ we proceed as follows. Since $b\in BMO(\mathbb R^n)$, as before, a straightforward computation shows that
\begin{equation}\label{BMO2}
\big|b_{2^{j+1}B}-b_{B}\big|\leq C\cdot(j+1)\|b\|_*.
\end{equation}
Thus, by the inequalities (\ref{BMO2}), (\ref{wq(B)/B}), the $A_1$ condition and the fact that $t\leq\Phi(t)$, we obtain
\begin{equation*}
\begin{split}
I'_6&\leq C\cdot w^q(B)^{1/q-\kappa}\sum_{j=1}^\infty(j+1)\|b\|_*\cdot\frac{1}{|2^{j+1}B|^{1-\alpha/n}}\int_{2^{j+1}B}\frac{|f(y)|}{\sigma}\,dy\\
&= C\cdot\|b\|_*\sum_{j=1}^\infty(j+1)\frac{w^q(B)^{1/q-\kappa}}{w^q(2^{j+1}B)^{1/q-\kappa}}\cdot\frac{1}{w^q(2^{j+1}B)^\kappa}
\cdot\frac{w^q(2^{j+1}B)^{1/q}}{|2^{j+1}B|^{1/q}}
\int_{2^{j+1}B}\frac{|f(y)|}{\sigma}\,dy\\
&\leq C\cdot\|b\|_*\sum_{j=1}^\infty(j+1)\frac{w^q(B)^{1/q-\kappa}}{w^q(2^{j+1}B)^{1/q-\kappa}}\cdot\frac{1}{w^q(2^{j+1}B)^\kappa}
\cdot\frac{w(2^{j+1}B)}{|2^{j+1}B|}
\int_{2^{j+1}B}\frac{|f(y)|}{\sigma}\,dy\\
&\leq C\cdot\|b\|_*\sum_{j=1}^\infty(j+1)\frac{w^q(B)^{1/q-\kappa}}{w^q(2^{j+1}B)^{1/q-\kappa}}\cdot\frac{1}{w^q(2^{j+1}B)^\kappa}
\int_{2^{j+1}B}\frac{|f(y)|}{\sigma}\cdot w(y)\,dy\\
&\leq C\cdot\sup_B\left\{\frac{1}{w^q(B)^\kappa}
\int_{B}\Phi\left(\frac{|f(y)|}{\sigma}\right)\cdot w(y)\,dy\right\}\times\sum_{j=1}^\infty(j+1)\cdot\frac{w^q(B)^{1/q-\kappa}}{w^q(2^{j+1}B)^{1/q-\kappa}}.
\end{split}
\end{equation*}
Moreover, since $w^q\in A_1\subset A_\infty$, by using the inequality (\ref{compare}) again, we have
\begin{align}
\sum_{j=1}^\infty(j+1)\cdot\frac{w^q(B)^{1/q-\kappa}}{w^q(2^{j+1}B)^{1/q-\kappa}}
&\leq C\sum_{j=1}^\infty(j+1)\cdot\left(\frac{|B|}{|2^{j+1}B|}\right)^{\delta^*(1/q-\kappa)}\notag\\
&\leq C\sum_{j=1}^\infty(j+1)\cdot\left(\frac{1}{2^{(j+1)n}}\right)^{\delta^*(1/q-\kappa)}\leq C,
\end{align}
which in turn gives that
\begin{equation*}
I'_6\leq C\cdot\sup_B\left\{\frac{1}{w^q(B)^\kappa}
\int_{B}\Phi\left(\frac{|f(y)|}{\sigma}\right)\cdot w(y)\,dy\right\}.
\end{equation*}
Combining all the above estimates, we are done.
\end{proof}

\begin{proof}[Proof of Theorem $\ref{mainthm:4}$]
For any ball $B=B(x_0,r)\subseteq\mathbb R^n$ with $x_0\in\mathbb R^n$ and $r>0$, we set $f=f\cdot\chi_{_{2B}}+f\cdot\chi_{_{(2B)^c}}:=f_1+f_2$. Then for each fixed $\sigma>0$, we have
\begin{equation*}
\begin{split}
&\left(\frac{1}{\Theta^q(r)}\cdot\big|\big\{x\in B:\big|[b,\mathcal T_\alpha](f)(x)\big|>\sigma\big\}\big|\right)^{1/q}\\
\leq &\left(\frac{1}{\Theta^q(r)}\cdot\big|\big\{x\in B:\big|[b,\mathcal T_\alpha](f_1)(x)\big|>\sigma/2\big\}\big|\right)^{1/q}\\
&+\left(\frac{1}{\Theta^q(r)}\cdot\big|\big\{x\in B:\big|[b,\mathcal T_\alpha](f_2)(x)\big|>\sigma/2\big\}\big|\right)^{1/q}\\
:=&J'_1+J'_2.
\end{split}
\end{equation*}
We consider the term $J'_1$ first. The assumption (\ref{assump4}) and the inequality (\ref{doubling}) yield that
\begin{equation*}
\begin{split}
J'_1&\leq  C_0\cdot\frac{1}{\Theta(r)}\int_{B}\Phi\left(\frac{|f_1(x)|}{\sigma}\right)dx\\
&= C_0\cdot\frac{1}{\Theta(r)}\int_{2B}\Phi\left(\frac{|f(x)|}{\sigma}\right)dx\\
&= C_0\cdot\frac{\Theta(2r)}{\Theta(r)}\cdot\frac{1}{\Theta(2r)}
\int_{B(x_0,2r)}\Phi\left(\frac{|f(x)|}{\sigma}\right)dx\\
&\leq C\cdot\sup_{r>0;B(x_0,r)}\left\{\frac{1}{\Theta(r)}
\int_{B(x_0,r)}\Phi\left(\frac{|f(x)|}{\sigma}\right)dx\right\}.
\end{split}
\end{equation*}
We now turn our attention to the estimate of $J'_2$. Recalling that the following estimate holds for given $0<\alpha<n$ and any $x\in B$,
\begin{equation*}
\big|[b,\mathcal T_\alpha](f_2)(x)\big|\leq\widetilde\mu(x)+\widetilde\nu(x),
\end{equation*}
where
\begin{equation*}
\widetilde\mu(x)=c_3\big|b(x)-b_{B}\big|\cdot\int_{\mathbb R^n}\frac{|f_2(y)|}{|x-y|^{n-\alpha}}dy,
\end{equation*}
and
\begin{equation*}
\widetilde\nu(x)=c_3\int_{\mathbb R^n}\frac{|b(y)-b_{B}|\cdot|f_2(y)|}{|x-y|^{n-\alpha}}dy.
\end{equation*}
Thus, we have
\begin{equation*}
\begin{split}
J'_2\leq&\left(\frac{1}{\Theta^q(r)}\cdot \big|\big\{x\in B:\widetilde\mu(x)>\sigma/4\big\}\big|\right)^{1/q}
+\left(\frac{1}{\Theta^q(r)}\cdot \big|\big\{x\in B:\widetilde\nu(x)>\sigma/4\big\}\big|\right)^{1/q}\\
:=&J'_3+J'_4.
\end{split}
\end{equation*}
Using the previous pointwise estimate (\ref{Ta(f2)}), Chebyshev's inequality, John--Nirenberg's inequality and the fact that $1/q=1-\alpha/n$, we conclude that
\begin{equation*}
\begin{split}
J'_3&\leq\frac{1}{\Theta(r)}\cdot\frac{\,4\,}{\sigma}\left(\int_B \widetilde\mu(x)^q\,dx\right)^{1/q}\\
&\leq \frac{C}{\Theta(r)}\sum_{j=1}^\infty\frac{1}{|2^{j+1}B|^{1-\alpha/n}}\int_{2^{j+1}B}\frac{|f(y)|}{\sigma}\,dy
\times\left\{|B|\cdot\frac{1}{|B|}\int_B\big|b(x)-b_{B}\big|^qdx\right\}^{1/q}\\
&\leq C\|b\|_*\sum_{j=1}^\infty\frac{|B|^{1/q}}{|2^{j+1}B|^{1/q}}\cdot\frac{\Theta(2^{j+1}r)}{\Theta(r)}\cdot
\frac{1}{\Theta(2^{j+1}r)}\int_{B(x_0,2^{j+1}r)}\frac{|f(y)|}{\sigma}\,dy.
\end{split}
\end{equation*}
Observe that $t\leq t\cdot(1+\log^+t)=\Phi(t)$, we get
\begin{equation*}
\begin{split}
J'_3&\leq C\|b\|_*\sum_{j=1}^\infty\frac{|B|^{1/q}}{|2^{j+1}B|^{1/q}}\cdot\frac{\Theta(2^{j+1}r)}{\Theta(r)}
\times\frac{1}{\Theta(2^{j+1}r)}\int_{B(x_0,2^{j+1}r)}\Phi\left(\frac{|f(y)|}{\sigma}\right)dy\\
&\leq C\cdot\sup_{r>0;B(x_0,r)}\left\{\frac{1}{\Theta(r)}
\int_{B(x_0,r)}\Phi\left(\frac{|f(y)|}{\sigma}\right)dy\right\}\times\sum_{j=1}^\infty\frac{|B|^{1/q}}{|2^{j+1}B|^{1/q}}\cdot\frac{\Theta(2^{j+1}r)}{\Theta(r)}.
\end{split}
\end{equation*}
Noting that $0<D(\Theta)<2^{n/q}$, then by using the doubling condition (\ref{doubling}) of $\Theta$, we are able to verify that
\begin{equation}\label{theta2<C}
\sum_{j=1}^\infty\frac{|B|^{1/q}}{|2^{j+1}B|^{1/q}}\cdot\frac{\Theta(2^{j+1}r)}{\Theta(r)}\leq
\sum_{j=1}^\infty\left(\frac{D(\Theta)}{2^{n/q}}\right)^{j+1}\leq C.
\end{equation}
Hence
\begin{equation*}
J'_3\leq C\cdot\sup_{r>0;B(x_0,r)}\left\{\frac{1}{\Theta(r)}
\int_{B(x_0,r)}\Phi\left(\frac{|f(y)|}{\sigma}\right)dy\right\}.
\end{equation*}
Applying the previous pointwise estimate (\ref{[b,Ta](f2)}) and Chebyshev's inequality, we have
\begin{equation*}
\begin{split}
J'_4&\leq\frac{1}{\Theta(r)}\cdot\frac{\,4\,}{\sigma}\left(\int_B \widetilde\nu(x)^q\,dx\right)^{1/q}\\
&\leq\frac{|B|^{1/q}}{\Theta(r)}\cdot\frac{C}{\sigma}\sum_{j=1}^\infty\frac{1}{|2^{j+1}B|^{1-\alpha/n}}\int_{2^{j+1}B}
\big|b(y)-b_{B}\big|\cdot\big|f(y)\big|\,dy\\
&\leq\frac{|B|^{1/q}}{\Theta(r)}\cdot\frac{C}{\sigma}\sum_{j=1}^\infty\frac{1}{|2^{j+1}B|^{1-\alpha/n}}\int_{2^{j+1}B}
\big|b(y)-b_{2^{j+1}B}\big|\cdot\big|f(y)\big|\,dy\\
&+\frac{|B|^{1/q}}{\Theta(r)}\cdot\frac{C}{\sigma}\sum_{j=1}^\infty\frac{1}{|2^{j+1}B|^{1-\alpha/n}}\int_{2^{j+1}B}
\big|b_{2^{j+1}B}-b_B\big|\cdot\big|f(y)\big|\,dy\\
&:=J'_5+J'_6.
\end{split}
\end{equation*}
For the term $J'_5$, notice that the inequality $\Phi(a\cdot b)\leq\Phi(a)\cdot\Phi(b)$ holds for any $a,b>0$, when $\Phi(t)=t\cdot(1+\log^+t)$. We then use the generalized H\"older's inequality (\ref{holder}), (\ref{exp}) and (\ref{equiv norm}) together with (\ref{theta2<C}) to obtain
\begin{equation*}
\begin{split}
J'_5&\leq \frac{|B|^{1/q}}{\Theta(r)}\cdot\frac{C}{\sigma}\sum_{j=1}^\infty\big|2^{j+1}B\big|^{\alpha/n}\cdot\big\|b-b_{2^{j+1}B}\big\|_{\exp L,2^{j+1}B}\big\|f\big\|_{L\log L,2^{j+1}B}\\
&\leq \frac{C\|b\|_*}{\sigma}\cdot\frac{|B|^{1/q}}{\Theta(r)}\sum_{j=1}^\infty\big|2^{j+1}B\big|^{\alpha/n}\times
\inf_{\eta>0}\left\{\eta+\frac{\eta}{|2^{j+1}B|}\int_{2^{j+1}B}\Phi\left(\frac{|f(z)|}{\eta}\right)dz\right\}\\
&\leq \frac{C\|b\|_*}{\sigma}\cdot \frac{|B|^{1/q}}{\Theta(r)}\sum_{j=1}^\infty\big|2^{j+1}B\big|^{\alpha/n}\times
\left\{\frac{\sigma\cdot\Theta(2^{j+1}r)}{|2^{j+1}B|}+\frac{\sigma}{|2^{j+1}B|}\int_{2^{j+1}B}\Phi\left(\frac{|f(z)|}{\sigma}\right)dz\right\}\\
&\leq C\|b\|_*\cdot\left[1+\sup_{r>0;B(x_0,r)}\left\{\frac{1}{\Theta(r)}
\int_{B(x_0,r)}\Phi\left(\frac{|f(z)|}{\sigma}\right)dz\right\}\right]\\
&\times\sum_{j=1}^\infty\frac{|B|^{1/q}}{|2^{j+1}B|^{1/q}}\cdot\frac{\Theta(2^{j+1}r)}{\Theta(r)}\\
&\leq C\cdot\sup_{r>0;B(x_0,r)}\left\{\frac{1}{\Theta(r)}
\int_{B(x_0,r)}\Phi\left(\frac{|f(z)|}{\sigma}\right)dz\right\}.
\end{split}
\end{equation*}
For the last term $J'_6$, in view of the inequality (\ref{BMO2}) and the fact that $t\leq\Phi(t)$, we get
\begin{equation*}
\begin{split}
J'_6&\leq C\cdot\frac{|B|^{1/q}}{\Theta(r)}\sum_{j=1}^\infty(j+1)\|b\|_*\cdot\frac{1}{|2^{j+1}B|^{1-\alpha/n}}\int_{2^{j+1}B}\frac{|f(y)|}{\sigma}\,dy\\
&\leq C\cdot\frac{|B|^{1/q}}{\Theta(r)}\sum_{j=1}^\infty(j+1)\|b\|_*\cdot\frac{\Theta(2^{j+1}r)}{|2^{j+1}B|^{1-\alpha/n}}\cdot
\frac{1}{\Theta(2^{j+1}r)}\int_{B(x_0,2^{j+1}r)}\Phi\left(\frac{|f(y)|}{\sigma}\right)dy\\
&\leq C\cdot\sup_{r>0;B(x_0,r)}\left\{\frac{1}{\Theta(r)}
\int_{B(x_0,r)}\Phi\left(\frac{|f(y)|}{\sigma}\right)dy\right\}\\
&\times\sum_{j=1}^\infty(j+1)\cdot\frac{|B|^{1/q}}{|2^{j+1}B|^{1/q}}\cdot\frac{\Theta(2^{j+1}r)}{\Theta(r)}.
\end{split}
\end{equation*}
Moreover, by using the doubling condition (\ref{doubling}) of $\Theta$ again and the fact that $0<D(\Theta)<2^{n/q}$, we find that
\begin{equation}\label{final}
\sum_{j=1}^\infty(j+1)\cdot\frac{|B|^{1/q}}{|2^{j+1}B|^{1/q}}\cdot\frac{\Theta(2^{j+1}r)}{\Theta(r)}\leq C\sum_{j=1}^\infty(j+1)\cdot\left(\frac{D(\Theta)}{2^{n/q}}\right)^{j+1}\leq C.
\end{equation}
Substituting the above inequality (\ref{final}) into the term $J'_6$, we finally obtain
\begin{equation*}
J'_6\leq C\cdot\sup_{r>0;B(x_0,r)}\left\{\frac{1}{\Theta(r)}
\int_{B(x_0,r)}\Phi\left(\frac{|f(y)|}{\sigma}\right)dy\right\}.
\end{equation*}
Summing up all the above estimates, we finish the proof of the main theorem.
\end{proof}

\section{Some applications}

In this section, we will give some applications of our main theorems to several integral operators such as $\theta$-type Calder\'on--Zygmund operators, Marcinkiewicz integral operators, Littlewood--Paley operators, Bochner--Riesz means, fractional maximal functions and fractional integrals.

\subsection{$\theta$-type Calder\'on--Zygmund operators}

Calder\'on--Zygmund singular integral operators and their generalizations on the Euclidean space $\mathbb R^n$ have been extensively studied (see \cite{duoand,garcia,stein2,yabuta} for instance). In particular, Yabuta \cite{yabuta} introduced certain $\theta$-type Calder\'on--Zygmund operators to facilitate his study of certain classes of pseudo-differential operator. Let $\theta$ be a non-negative, non-decreasing function on $(0,+\infty)$ with
\begin{equation*}
\int_0^1\frac{\theta(t)\cdot|\log t|}{t}dt<\infty.
\end{equation*}
A measurable function $K$ on $\mathbb R^n\times\mathbb R^n\backslash\{(x,x):x\in\mathbb R^n\}$ is said to be a $\theta$-type kernel if it satisfies

(i) $|K(x,y)|\leq C\cdot|x-y|^{-n}$, for any $x\neq y$;

(ii) $|K(x,y)-K(z,y)|+|K(y,x)-K(y,z)|\leq C\cdot\theta({|x-z|}/{|x-y|})|x-y|^{-n}$, for $|x-z|<|x-y|/2$.

Let $T_\theta$ be a linear operator from $\mathscr S(\mathbb R^n)$ into its dual $\mathscr S'(\mathbb R^n)$. We say that $T_\theta$ is a $\theta$-type Calder\'on--Zygmund operator if

(1) $T_\theta$ can be extended to be a bounded operator on $L^2(\mathbb R^n)$;

(2) There is a $\theta$-type kernel $K$ such that $Tf(x)=\int_{\mathbb R^n}K(x,y)f(y)\,dy$ for all $f\in C^\infty_0(\mathbb R^n)$ and for all $x\notin supp\, f$, where $C^\infty_0(\mathbb R^n)$ is the space consisting of all infinitely differentiable functions on $\mathbb R^n$ with compact supports. If $b\in BMO(\mathbb R^n)$, we define the commutator $[b,T_\theta]$ to be the operator
$$[b,T_\theta]f(x)=b(x)\cdot T_\theta f(x)-T_\theta(b\cdot f)(x)=\int_{\mathbb R^n}[b(x)-b(y)]K(x,y)f(y)\,dy.$$
The following endpoint estimates for commutator of the $\theta$-type Calder\'on--Zygmund operator were established in \cite{liu,zhang2}.

\begin{theorem}[\cite{zhang2}]
Let $w\in A_1$ and $b\in BMO(\mathbb R^n)$. Then for all $\sigma>0$, there is a constant $C_0>0$ independent of $f$ and $\sigma$ such that
\begin{equation*}
w\big(\big\{x\in\mathbb R^n:\big|[b,T_\theta](f)(x)\big|>\sigma\big\}\big)\leq C_0\int_{\mathbb R^n}\Phi\left(\frac{|f(x)|}{\sigma}\right)\cdot w(x)\,dx,
\end{equation*}
where $\Phi(t)=t(1+\log^+t)$.
\end{theorem}

\begin{theorem}[\cite{liu}]
Let $b\in BMO(\mathbb R^n)$. Then for all $\sigma>0$, there is a constant $C_0>0$ independent of $f$ and $\sigma$ such that
\begin{equation*}
\big|\big\{x\in\mathbb R^n:\big|[b,T_\theta](f)(x)\big|>\sigma\big\}\big|\leq C_0\int_{\mathbb R^n}\Phi\left(\frac{|f(x)|}{\sigma}\right)dx,
\end{equation*}
where $\Phi(t)=t(1+\log^+t)$.
\end{theorem}

Then, from Theorem $\ref{mainthm:1}$ and Theorem $\ref{mainthm:2}$, we immediately get the following:
\begin{corollary}
Let $0<\kappa<1$, $w\in A_1$ and $b\in BMO(\mathbb R^n)$. Then for any given $\sigma>0$ and any ball $B$, there exists a constant $C>0$ independent of $f$, $B$ and $\sigma$ such that
\begin{equation*}
\frac{1}{w(B)^\kappa}\cdot w\big(\big\{x\in B:\big|[b,T_\theta](f)(x)\big|>\sigma\big\}\big)\leq C\cdot\sup_B\frac{1}{w(B)^\kappa}
\int_{B}\Phi\left(\frac{|f(x)|}{\sigma}\right)\cdot w(x)\,dx,
\end{equation*}
where $\Phi(t)=t(1+\log^+t)$.
\end{corollary}

\begin{corollary}
Let $b\in BMO(\mathbb R^n)$. Suppose that $\Theta$ satisfies $(\ref{doubling})$ and $0<D(\Theta)<2^n$, then for any given $\sigma>0$ and any ball $B(x_0,r)$, there exists a constant $C>0$ independent of $f$, $B(x_0,r)$ and $\sigma$ such that
\begin{equation*}
\frac{1}{\Theta(r)}\cdot\big|\big\{x\in B(x_0,r):\big|[b,T_\theta](f)(x)\big|>\sigma\big\}\big|\leq C\cdot\sup_{r>0}\frac{1}{\Theta(r)}
\int_{B(x_0,r)}\Phi\left(\frac{|f(x)|}{\sigma}\right)dx,
\end{equation*}
where $\Phi(t)=t(1+\log^+t)$.
\end{corollary}

\subsection{Marcinkiewicz integral operators}

Suppose that $S^{n-1}$ is the unit sphere in $\mathbb R^n$($n\ge2$) equipped with the normalized Lebesgue measure $d\sigma$. Let $\Omega$ be a homogeneous function of degree zero on $\mathbb R^n$ satisfying $\Omega\in L^1(S^{n-1})$ and $\int_{S^{n-1}}\Omega(x')\,d\sigma(x')=0$, where $x'=x/{|x|}$ for any $x\neq0$. Then the Marcinkiewicz integral of higher dimension is defined by
\begin{equation*}
\mu_{\Omega}(f)(x)=\left(\int_0^\infty\big|F_{\Omega,t}(x)\big|^2\frac{dt}{t^{3}}\right)^{1/2},
\end{equation*}
where
\begin{equation*}
F_{\Omega,t}(x)=\int_{|x-y|\le t}\frac{\Omega(x-y)}{|x-y|^{n-1}}f(y)\,dy.
\end{equation*}
For $b\in BMO(\mathbb R^n)$, the commutator operator $[b,\mu_{\Omega}]$ is defined by (see \cite{ding2})
\begin{equation*}
[b,\mu_{\Omega}](f)(x)=\left(\int_0^\infty\big|F^b_{\Omega,t}(x)\big|^2\frac{dt}{t^{3}}\right)^{1/2},
\end{equation*}
where
\begin{equation*}
F^b_{\Omega,t}(x)=\int_{|x-y|\le t}\frac{\Omega(x-y)}{|x-y|^{n-1}}[b(x)-b(y)]f(y)\,dy.
\end{equation*}
For $0<\alpha\leq1$, we say that $\Omega\in Lip_\alpha(S^{n-1})$, if there exists a constant $L>0$ such that
\begin{equation*}
\big|\Omega(x')-\Omega(y')\big|\leq L|x'-y'|^\alpha, \quad \mbox{for any }\, x',y'\in S^{n-1}.
\end{equation*}
Let $\mathbb H$ be the Banach space
\begin{equation*}
\mathbb H=\left\{h:\|h\|=\left(\int_0^\infty|h(t)|^2\frac{dt}{t^3}\right)^{1/2}<\infty\right\}.
\end{equation*}
Then, it is clear that $[b,\mu_{\Omega}](f)(x)=\|F^b_{\Omega,t}(x)\|$. By Minkowski's inequality and the condition on $\Omega$, we can get
\begin{equation*}
\begin{split}
\big|[b,\mu_{\Omega}](f)(x)\big|&\leq\int_{\mathbb R^n}\frac{|\Omega(x-y)|}{|x-y|^{n-1}}\big|b(x)-b(y)\big|\cdot\big|f(y)\big|\left(\int_{|x-y|}^\infty\frac{dt}{t^3}\right)^{1/2}dy\\
&\leq c_2\int_{\mathbb R^n}\frac{|b(x)-b(y)|\cdot|f(y)|}{|x-y|^n}dy,
\end{split}
\end{equation*}
where $c_2$ is an absolute constant independent of $f$ and $x\in\mathbb R^n$. Thus, $[b,\mu_{\Omega}]$ satisfies the condition (\ref{sublinear commutator}). Moreover, in \cite{ding2}, Ding et al. considered the weighted weak $L\log L$-type estimate for the commutator $[b,\mu_{\Omega}]$ and proved:

\begin{theorem}[\cite{ding2}]
Let $0<\alpha\leq1$, $\Omega\in Lip_\alpha(S^{n-1})$, $w\in A_1$ and $b\in BMO(\mathbb R^n)$. Then for all $\sigma>0$, there is a constant $C_0>0$ independent of $f$ and $\sigma$ such that
\begin{equation*}
w\big(\big\{x\in\mathbb R^n:\big|[b,\mu_{\Omega}](f)(x)\big|>\sigma\big\}\big)\leq C_0\int_{\mathbb R^n}\Phi\left(\frac{|f(x)|}{\sigma}\right)\cdot w(x)\,dx,
\end{equation*}
where $\Phi(t)=t(1+\log^+t)$.
\end{theorem}

In particular, we have the following estimate if $w$ is taken to be a constant function.

\begin{theorem}[\cite{ding2}]
Let $0<\alpha\leq1$, $\Omega\in Lip_\alpha(S^{n-1})$ and $b\in BMO(\mathbb R^n)$. Then for all $\sigma>0$, there is a constant $C_0>0$ independent of $f$ and $\sigma$ such that
\begin{equation*}
\big|\big\{x\in\mathbb R^n:\big|[b,\mu_{\Omega}](f)(x)\big|>\sigma\big\}\big|\leq C_0\int_{\mathbb R^n}\Phi\left(\frac{|f(x)|}{\sigma}\right)dx,
\end{equation*}
where $\Phi(t)=t(1+\log^+t)$.
\end{theorem}

As a consequence of Theorem $\ref{mainthm:1}$ and Theorem $\ref{mainthm:2}$, we obtain the following results:
\begin{corollary}
Let $0<\alpha\leq1$, $\Omega\in Lip_\alpha(S^{n-1})$, $0<\kappa<1$, $w\in A_1$ and $b\in BMO(\mathbb R^n)$. Then for any given $\sigma>0$ and any ball $B$, there exists a constant $C>0$ independent of $f$, $B$ and $\sigma$ such that
\begin{equation*}
\frac{1}{w(B)^\kappa}\cdot w\big(\big\{x\in B:\big|[b,\mu_{\Omega}](f)(x)\big|>\sigma\big\}\big)\leq C\cdot\sup_B\frac{1}{w(B)^\kappa}
\int_{B}\Phi\left(\frac{|f(x)|}{\sigma}\right)\cdot w(x)\,dx,
\end{equation*}
where $\Phi(t)=t(1+\log^+t)$.
\end{corollary}

\begin{corollary}
Let $0<\alpha\leq1$, $\Omega\in Lip_\alpha(S^{n-1})$ and $b\in BMO(\mathbb R^n)$. Suppose that $\Theta$ satisfies $(\ref{doubling})$ and $0<D(\Theta)<2^n$, then for any given $\sigma>0$ and any ball $B(x_0,r)$, there exists a constant $C>0$ independent of $f$, $B(x_0,r)$ and $\sigma$ such that
\begin{equation*}
\frac{1}{\Theta(r)}\cdot\big|\big\{x\in B(x_0,r):\big|[b,\mu_{\Omega}](f)(x)\big|>\sigma\big\}\big|\leq C\cdot\sup_{r>0}\frac{1}{\Theta(r)}
\int_{B(x_0,r)}\Phi\left(\frac{|f(x)|}{\sigma}\right)dx,
\end{equation*}
where $\Phi(t)=t(1+\log^+t)$.
\end{corollary}

\subsection{Littlewood--Paley operators}

Let $\varepsilon>0$ and $\psi$ be a fixed function which satisfies the following properties:

(1) $\psi\in L^1(\mathbb R^n)$ and $\int_{\mathbb R^n}\psi(x)\,dx=0$;

(2) $\psi(x)\leq C\cdot(1+|x|)^{-(n+1)}$;

(3) $\big|\psi(x+y)-\psi(x)\big|\leq C\cdot|y|^{\varepsilon}(1+|x|)^{-(n+1+\varepsilon)}$ when $2|y|<|x|$.

We set $\psi_t(x)=t^{-n}\psi(x/t)$ and $\Gamma(x)=\big\{(y,t)\in{\mathbb R}^{n+1}_+:|x-y|<t\big\}$. The Littlewood--Paley $g$-function, Lusin area integrals and the $g^*_\lambda$-function will be defined respectively by (see \cite{torchinsky})
\begin{equation*}
g_\psi(f)(x)=\left(\int_0^\infty\big|\psi_t*f(x)\big|^2\frac{dt}{t}\right)^{1/2},
\end{equation*}
\begin{equation*}
S_{\psi}(f)(x)=\left(\iint_{\Gamma(x)}\big|\psi_t*f(y)\big|^2\frac{dydt}{t^{n+1}}\right)^{1/2},
\end{equation*}
and
\begin{equation*}
g^*_{\lambda,\psi}(f)(x)=\left(\iint_{{\mathbb R}^{n+1}_+}\left(\frac t{t+|x-y|}\right)^{\lambda n}\big|\psi_t*f(y)\big|^2\frac{dydt}{t^{n+1}}\right)^{1/2}, \quad \lambda>1.
\end{equation*}
For $b\in BMO(\mathbb R^n)$, we will consider the commutators generated by $b$ and Littlewood--Paley operators, which are defined respectively by the following expressions (see \cite{xue}):
\begin{equation*}
\big[b,g_\psi\big](f)(x)=\left(\int_0^\infty\bigg|\int_{\mathbb R^n}\big[b(x)-b(y)\big]\psi_t(x-y)f(y)\,dy\bigg|^2\frac{dt}{t}\right)^{1/2},
\end{equation*}
\begin{equation*}
\big[b,S_\psi\big](f)(x)=\left(\iint_{\Gamma(x)}\bigg|\int_{\mathbb R^n}\big[b(x)-b(z)\big]\psi_t(y-z)f(z)\,dz\bigg|^2\frac{dydt}{t^{n+1}}\right)^{1/2},
\end{equation*}
and
\begin{equation*}
\begin{split}
&\big[b,g^*_{\lambda,\psi}\big](f)(x)\\
=&\left(\iint_{{\mathbb R}^{n+1}_+}\left(\frac t{t+|x-y|}\right)^{\lambda n}\bigg|\int_{\mathbb R^n}\big[b(x)-b(z)\big]\psi_t(y-z)f(z)\,dz\bigg|^2\frac{dydt}{t^{n+1}}\right)^{1/2}, \lambda>1.
\end{split}
\end{equation*}
Let $\mathbb H$ be the Banach space
\begin{equation*}
\mathbb H=\left\{h:\|h\|=\left(\int_0^\infty|h(t)|^2\frac{dt}{t}\right)^{1/2}<\infty\right\}
\end{equation*}
or
\begin{equation*}
\mathbb H=\left\{h:\|h\|=\left(\iint_{{\mathbb R}^{n+1}_+}|h(y,t)|^2\frac{dydt}{t^{n+1}}\right)^{1/2}<\infty\right\}.
\end{equation*}
If we set
\begin{equation*}
F^b_{\psi,t}(x)=\int_{\mathbb R^n}\big[b(x)-b(y)\big]\psi_t(x-y)f(y)\,dy,
\end{equation*}
\begin{equation*}
F^b_{\psi,t}(x,y)=\int_{\mathbb R^n}\big[b(x)-b(z)\big]\psi_t(y-z)f(z)\,dz,
\end{equation*}
and denote the characteristic function of $\Gamma(x)$ by $\chi_{\Gamma(x)}$, then, for each fixed $x\in\mathbb R^n$, it is easy to see that $$\big[b,g_\psi\big](f)(x)=\big\|F^b_{\psi,t}(x)\big\|,\quad\big[b,S_\psi\big](f)(x)=\big\|\chi_{\Gamma(x)}\cdot F^b_{\psi,t}(x,y)\big\|,$$
and
$$\big[b,g^*_{\lambda,\psi}\big](f)(x)=\left\|\left(\frac t{t+|x-y|}\right)^{\lambda n/2}\cdot F^b_{\psi,t}(x,y)\right\|.$$
By using Minkowski's inequality and the condition on $\psi$, we can get
\begin{equation*}
\begin{split}
\big|\big[b,g_\psi\big](f)(x)\big|&\leq c_2\int_{\mathbb R^n}\big|b(x)-b(y)\big|\cdot\big|f(y)\big|\left(\int_{0}^\infty\left(\frac{1}{t^n}\cdot\frac{1}{[1+t^{-1}|x-y|]^{n+1}}\right)^2\frac{dt}{t}\right)^{1/2}dy\\
&\leq c_2\int_{\mathbb R^n}\frac{|b(x)-b(y)|\cdot|f(y)|}{|x-y|^n}dy,
\end{split}
\end{equation*}
Similarly, we can also prove
$$\big|\big[b,S_\psi\big](f)(x)\big|\leq c_2\int_{\mathbb R^n}\frac{|b(x)-b(y)|\cdot|f(y)|}{|x-y|^n}dy,$$
and
$$\big|\big[b,g^*_{\lambda,\psi}\big](f)(x)\big|\leq c_2\int_{\mathbb R^n}\frac{|b(x)-b(y)|\cdot|f(y)|}{|x-y|^n}dy,$$
where $c_2$ is an absolute constant independent of $f$ and $x\in\mathbb R^n$. Thus, $[b,g_{\psi}]$, $[b,S_{\psi}]$ and $[b,g^*_{\lambda,\psi}]$ all satisfy the condition (\ref{sublinear commutator}). For the endpoint estimates for these commutator operators $[b,g_{\psi}]$, $[b,S_{\psi}]$ and $[b,g^*_{\lambda,\psi}]$ in the weighted Lebesgue space $L^1_w(\mathbb R^n)$,
when $b\in BMO(\mathbb R^n)$ and $w\in A_1$, it was proved by Xue and Ding in \cite{xue} that

\begin{theorem}[\cite{xue}]
Let $\lambda>3$, $w\in A_1$, $b\in BMO(\mathbb R^n)$ and $\psi$ be a function on $\mathbb R^n$ satisfying $(1)-(3)$ mentioned above. Then for all $\sigma>0$, there is a constant $C_0>0$ independent of $f$ and $\sigma$ such that
\begin{equation*}
w\big(\big\{x\in\mathbb R^n:\big|[b,T_\psi](f)(x)\big|>\sigma\big\}\big)\leq C_0\int_{\mathbb R^n}\Phi\left(\frac{|f(x)|}{\sigma}\right)\cdot w(x)\,dx,
\end{equation*}
where $\Phi(t)=t(1+\log^+t)$ and $T_\psi$ is $g_\psi$ or $S_\psi$ or $g^*_{\lambda,\psi}$.
\end{theorem}

\begin{theorem}[\cite{xue}]
Let $\lambda>3$, $b\in BMO(\mathbb R^n)$ and $\psi$ be a function on $\mathbb R^n$ satisfying $(1)-(3)$ mentioned above. Then for all $\sigma>0$, there is a constant $C_0>0$ independent of $f$ and $\sigma$ such that
\begin{equation*}
\big|\big\{x\in\mathbb R^n:\big|[b,T_\psi](f)(x)\big|>\sigma\big\}\big|\leq C_0\int_{\mathbb R^n}\Phi\left(\frac{|f(x)|}{\sigma}\right)dx,
\end{equation*}
where $\Phi(t)=t(1+\log^+t)$ and $T_\psi$ is $g_\psi$ or $S_\psi$ or $g^*_{\lambda,\psi}$.
\end{theorem}
Then, from Theorem $\ref{mainthm:1}$ and Theorem $\ref{mainthm:2}$, we will show that:

\begin{corollary}
Let $\lambda>3$, $0<\kappa<1$, $w\in A_1$, $b\in BMO(\mathbb R^n)$ and $\psi$ be a function on $\mathbb R^n$ satisfying $(1)-(3)$ mentioned above. Then for any given $\sigma>0$ and any ball $B$, there exists a constant $C>0$ independent of $f$, $B$ and $\sigma$ such that
\begin{equation*}
\frac{1}{w(B)^\kappa}\cdot w\big(\big\{x\in B:\big|[b,T_\psi](f)(x)\big|>\sigma\big\}\big)\leq C\cdot\sup_B\frac{1}{w(B)^\kappa}
\int_{B}\Phi\left(\frac{|f(x)|}{\sigma}\right)\cdot w(x)\,dx,
\end{equation*}
where $\Phi(t)=t(1+\log^+t)$ and $T_\psi$ is $g_\psi$ or $S_\psi$ or $g^*_{\lambda,\psi}$.
\end{corollary}

\begin{corollary}
Let $\lambda>3$, $b\in BMO(\mathbb R^n)$ and $\psi$ be a function on $\mathbb R^n$ satisfying $(1)-(3)$ mentioned above. Suppose that $\Theta$ satisfies $(\ref{doubling})$ and $0<D(\Theta)<2^n$, then for any given $\sigma>0$ and any ball $B(x_0,r)$, there exists a constant $C>0$ independent of $f$, $B(x_0,r)$ and $\sigma$ such that
\begin{equation*}
\frac{1}{\Theta(r)}\cdot\big|\big\{x\in B(x_0,r):\big|[b,T_\psi](f)(x)\big|>\sigma\big\}\big|\leq C\cdot\sup_{r>0}\frac{1}{\Theta(r)}
\int_{B(x_0,r)}\Phi\left(\frac{|f(x)|}{\sigma}\right)dx,
\end{equation*}
where $\Phi(t)=t(1+\log^+t)$ and $T_\psi$ is $g_\psi$ or $S_\psi$ or $g^*_{\lambda,\psi}$.
\end{corollary}

\subsection{Bochner--Riesz means}

The Bochner--Riesz means of order $\delta>0$ in $\mathbb R^n$ are defined initially for Schwartz functions in terms of Fourier transforms by
\begin{equation*}
\big(\widehat{T^\delta_R f}\big)(\xi)=\Big(1-\frac{|\xi|^2}{R^2}\Big)^\delta_+\widehat{f}(\xi), \quad 0<R<\infty,
\end{equation*}
where $\widehat{f}$ denotes the Fourier transform of $f$. We recall that the Bochner--Riesz means can be expressed as convolution operators (see \cite{lu1,stein1})
\begin{equation*}
T^\delta_Rf(x)=(\phi_{1/R}*f)(x),
\end{equation*}
where $\phi(x)=[(1-|\cdot|^2)^\delta_+]\mbox{\textasciicircum}(x)$ and $\phi_{1/R}(x)=R^n\cdot\phi(Rx)$. It is well known that the kernel $\phi$ can be represented as (see \cite{lu1,stein1})
\begin{equation*}
\phi(x)=\pi^{-\delta}\Gamma(\delta+1)|x|^{-(\frac n2+\delta)}J_{\frac n2+\delta}(2\pi|x|),
\end{equation*}
where $J_\mu(t)$ is the Bessel function
\begin{equation*}
J_\mu(t)=\frac{(\frac{t}{2})^\mu}{\Gamma(\mu+\frac12)\Gamma(\frac12)}\int_{-1}^1e^{its}(1-s^2)^{\mu-\frac12}\,ds.
\end{equation*}
Let $b\in BMO(\mathbb R^n)$ and $0<R<\infty$. Consider the commutator $\big[b,T^\delta_R\big]$ defined by
$$\big[b,T^\delta_R\big](f)(x)=b(x)\cdot T^\delta_Rf(x)-T^\delta_R(b\cdot f)(x)=\int_{\mathbb R^n}[b(x)-b(y)]\phi_{1/R}(x-y)f(y)\,dy.$$
The maximal operator $\big[b,T^\delta_*\big]$ associated with the commutator is defined by
\begin{equation*}
\big[b,T^\delta_*\big](f)(x)=\sup_{R>0}\Big|\big[b,T^\delta_R\big](f)(x)\Big|.
\end{equation*}
Let $\mathbb H$ be the space
\begin{equation*}
\mathbb H=\left\{h:\|h\|=\sup_{R>0}\big|h(R)\big|<\infty\right\}.
\end{equation*}
Then, it is clear that $\big[b,T^\delta_*\big](f)(x)=\|b(x)\cdot T^\delta_Rf(x)-T^\delta_R(b\cdot f)(x)\|$. If $\delta\geq{(n-1)}/2$, by the kernel estimates of $T^\delta_R$, we have
\begin{equation*}
\begin{split}
\big|\big[b,T^\delta_*\big](f)(x)\big|&\leq c_2\cdot\sup_{R>0}\int_{\mathbb R^n}\big|b(x)-b(y)\big|\frac{R^n}{(1+R|x-y|)^{\delta+\frac{n+1}{2}}}\cdot\big|f(y)\big|\,dy\\
&\leq c_2\int_{\mathbb R^n}\frac{|b(x)-b(y)|\cdot|f(y)|}{|x-y|^n}\,dy,
\end{split}
\end{equation*}
where $c_2$ is an absolute constant independent of $f$ and $x\in\mathbb R^n$. Thus, $\big[b,T^\delta_*\big]$ satisfies the condition (\ref{sublinear commutator}). Furthermore, in \cite{liulz}, Liu and Lu established weighted endpoint estimates of $L\log L$-type for maximal commutators of the Bochner--Riesz means.

\begin{theorem}[\cite{liulz}]
Let $\delta>{(n-1)}/2$, $w\in A_1$ and $b\in BMO(\mathbb R^n)$. Then for all $\sigma>0$, there is a constant $C_0>0$ independent of $f$ and $\sigma$ such that
\begin{equation*}
w\big(\big\{x\in\mathbb R^n:\big|\big[b,T^\delta_*\big](f)(x)\big|>\sigma\big\}\big)\leq C_0\int_{\mathbb R^n}\Phi\left(\frac{|f(x)|}{\sigma}\right)\cdot w(x)\,dx,
\end{equation*}
where $\Phi(t)=t(1+\log^+t)$.
\end{theorem}

\begin{theorem}[\cite{liulz}]
Let $\delta>{(n-1)}/2$ and $b\in BMO(\mathbb R^n)$. Then for all $\sigma>0$, there is a constant $C_0>0$ independent of $f$ and $\sigma$ such that
\begin{equation*}
\big|\big\{x\in\mathbb R^n:\big|\big[b,T^\delta_*\big](f)(x)\big|>\sigma\big\}\big|\leq C_0\int_{\mathbb R^n}\Phi\left(\frac{|f(x)|}{\sigma}\right)dx,
\end{equation*}
where $\Phi(t)=t(1+\log^+t)$.
\end{theorem}

As a consequence of Theorem $\ref{mainthm:1}$ and Theorem $\ref{mainthm:2}$, we can prove the following results:
\begin{corollary}
Let $\delta>{(n-1)}/2$, $0<\kappa<1$, $w\in A_1$ and $b\in BMO(\mathbb R^n)$. Then for any given $\sigma>0$ and any ball $B$, there exists a constant $C>0$ independent of $f$, $B$ and $\sigma$ such that
\begin{equation*}
\frac{1}{w(B)^\kappa}\cdot w\big(\big\{x\in B:\big|\big[b,T^\delta_*\big](f)(x)\big|>\sigma\big\}\big)\leq C\cdot\sup_B\frac{1}{w(B)^\kappa}
\int_{B}\Phi\left(\frac{|f(x)|}{\sigma}\right)\cdot w(x)\,dx,
\end{equation*}
where $\Phi(t)=t(1+\log^+t)$.
\end{corollary}

\begin{corollary}
Let $\delta>{(n-1)}/2$ and $b\in BMO(\mathbb R^n)$. Suppose that $\Theta$ satisfies $(\ref{doubling})$ and $0<D(\Theta)<2^n$, then for any given $\sigma>0$ and any ball $B(x_0,r)$, there exists a constant $C>0$ independent of $f$, $B(x_0,r)$ and $\sigma$ such that
\begin{equation*}
\frac{1}{\Theta(r)}\cdot\big|\big\{x\in B(x_0,r):\big|\big[b,T^\delta_*\big](f)(x)\big|>\sigma\big\}\big|\leq C\cdot\sup_{r>0}\frac{1}{\Theta(r)}
\int_{B(x_0,r)}\Phi\left(\frac{|f(x)|}{\sigma}\right)dx,
\end{equation*}
where $\Phi(t)=t(1+\log^+t)$.
\end{corollary}

\subsection{Fractional integrals}

For given $\alpha$, $0<\alpha<n$, the fractional integral operator (or the Riesz potential) $I_\alpha$ is defined by (see \cite{stein})
\begin{equation*}
I_\alpha f(x)=\frac{1}{\gamma(\alpha)}\int_{\mathbb R^n}\frac{f(y)}{|x-y|^{n-\alpha}}\,dy, \quad \mbox{ÆäÖÐ}\; \gamma(\alpha)=\frac{2^\alpha\pi^{\frac n2}\Gamma(\frac{\alpha}{2})}{\Gamma(\frac{n-\alpha}{2})}.
\end{equation*}
We also define the associated fractional maximal function with order $\alpha$ by
\begin{equation*}
M_{\alpha}(f)(x)=\sup_{x\in B}\frac{1}{|B|^{1-\frac{\alpha}{n}}}\int_B|f(y)|\,dy,
\end{equation*}
where the supremum is taken over all balls containing $x$. When $b\in BMO(\mathbb R^n)$, the commutators $[b,I_\alpha]$ and $[b,M_\alpha]$ are defined as
$$[b,I_\alpha]f(x)=b(x)\cdot I_\alpha f(x)-I_\alpha(b\cdot f)(x)=\int_{\mathbb R^n}[b(x)-b(y)]\cdot\frac{f(y)}{|x-y|^{n-\alpha}}\,dy,$$
$$[b,M_\alpha](f)(x)=\sup_{x\in B}\frac{1}{|B|^{1-\frac{\alpha}{n}}}\int_B\big|b(x)-b(y)\big|\cdot|f(y)|\,dy.$$
In \cite{cruz1,cruz2}, Cruz-Uribe and Fiorenza discussed the unweighted and weighted endpoint inequalities for commutators of fractional integrals and proved the following

\begin{theorem}
Let $0<\alpha<n$, $q=n/{(n-\alpha)}$, $w^q\in A_1$ and $b\in BMO(\mathbb R^n)$. Then for any given $\sigma>0$ and any bounded domain $\Omega\subset\mathbb R^n$, there is a constant $C_0>0$ which does not depend on $f$, $\Omega$ and $\sigma$ such that
\begin{equation*}
\left[w^q\big(\big\{x\in\Omega:\big|[b,I_\alpha](f)(x)\big|>\sigma\big\}\big)\right]^{1/q}\leq C_0\int_{\Omega}\Phi\left(\frac{|f(x)|}{\sigma}\right)\cdot w(x)\,dx,
\end{equation*}
where $\Phi(t)=t(1+\log^+t)$.
\end{theorem}

\begin{theorem}
Let $0<\alpha<n$, $q=n/{(n-\alpha)}$ and $b\in BMO(\mathbb R^n)$. Then for any given $\sigma>0$ and any bounded domain $\Omega\subset\mathbb R^n$, there is a constant $C_0>0$ which does not depend on $f$, $\Omega$ and $\sigma$ such that
\begin{equation*}
\big|\big\{x\in\Omega:\big|[b,I_\alpha](f)(x)\big|>\sigma\big\}\big|^{1/q}\leq C_0\int_{\Omega}\Phi\left(\frac{|f(x)|}{\sigma}\right)dx,
\end{equation*}
where $\Phi(t)=t(1+\log^+t)$.
\end{theorem}

Then, from Theorem $\ref{mainthm:3}$ and Theorem $\ref{mainthm:4}$, we immediately get the following:

\begin{corollary}
Let $0<\alpha<n$, $q=n/{(n-\alpha)}$, $0<\kappa<1/q$, $w^q\in A_1$ and $b\in BMO(\mathbb R^n)$. Then for any given $\sigma>0$ and any ball $B\subset\mathbb R^n$, there exists a constant $C>0$ independent of $f$, $B$ and $\sigma$ such that
\begin{equation*}
\begin{split}
&\left(\frac{1}{w^q(B)^{\kappa q}}\cdot w^q\big(\big\{x\in B:\big|[b,I_{\alpha}](f)(x)\big|>\sigma\big\}\big)\right)^{1/q}\\
\leq& C\cdot\sup_B\frac{1}{w^q(B)^\kappa}
\int_{B}\Phi\left(\frac{|f(x)|}{\sigma}\right)\cdot w(x)\,dx,
\end{split}
\end{equation*}
where $\Phi(t)=t(1+\log^+t)$.
\end{corollary}

\begin{corollary}
Let $0<\alpha<n$, $q=n/{(n-\alpha)}$ and $b\in BMO(\mathbb R^n)$. Suppose that $\Theta$ satisfies $(\ref{doubling})$ and $0<D(\Theta)<2^{n/q}$, then for any given $\sigma>0$ and any ball $B(x_0,r)\subset\mathbb R^n$, there exists a constant $C>0$ independent of $f$, $B(x_0,r)$ and $\sigma$ such that
\begin{equation*}
\begin{split}
&\left(\frac{1}{\Theta^q(r)}\cdot\big|\big\{x\in B(x_0,r):\big|[b,I_\alpha](f)(x)\big|>\sigma\big\}\big|\right)^{1/q}\\
\leq& C\cdot\sup_{r>0}\frac{1}{\Theta(r)}\int_{B(x_0,r)}\Phi\left(\frac{|f(x)|}{\sigma}\right)dx,
\end{split}
\end{equation*}
where $\Phi(t)=t(1+\log^+t)$.
\end{corollary}
It should be pointed out that $[b,M_{\alpha}](f)$ can be controlled pointwise
by $[b,I_{\alpha}](|f|)$ for any $f(x)$ (see \cite{ding1}). In fact, for any $0<\alpha<n$, $x\in\mathbb R^n$ and $r>0$, we have
\begin{equation*}
\begin{split}
[b,I_{\alpha}](|f|)(x)&\ge\int_{|y-x|\le r}\frac{|b(x)-b(y)|\cdot|f(y)|}{|x-y|^{n-\alpha}}\,dy\\
&\ge \frac{1}{r^{n-\alpha}}\int_{|y-x|\le r}|b(x)-b(y)|\cdot|f(y)|\,dy.
\end{split}
\end{equation*}
Taking the supremum for all $r>0$ on both sides of the above inequality, we get
\begin{equation*}
[b,M_{\alpha}](f)(x)\le[b,I_{\alpha}](|f|)(x),\quad\mbox{for all}\; x\in\mathbb R^n.
\end{equation*}

Hence, as a direct consequence of the above results, we have eventually obtained

\begin{corollary}
Let $0<\alpha<n$, $q=n/{(n-\alpha)}$, $0<\kappa<1/q$, $w^q\in A_1$ and $b\in BMO(\mathbb R^n)$. Then for any given $\sigma>0$ and any ball $B\subset\mathbb R^n$, there exists a constant $C>0$ independent of $f$, $B$ and $\sigma$ such that
\begin{equation*}
\begin{split}
&\left(\frac{1}{w^q(B)^{\kappa q}}\cdot w^q\big(\big\{x\in B:\big|[b,M_{\alpha}](f)(x)\big|>\sigma\big\}\big)\right)^{1/q}\\
\leq& C\cdot\sup_B\frac{1}{w^q(B)^\kappa}
\int_{B}\Phi\left(\frac{|f(x)|}{\sigma}\right)\cdot w(x)\,dx,
\end{split}
\end{equation*}
where $\Phi(t)=t(1+\log^+t)$.
\end{corollary}

\begin{corollary}
Let $0<\alpha<n$, $q=n/{(n-\alpha)}$ and $b\in BMO(\mathbb R^n)$. Suppose that $\Theta$ satisfies $(\ref{doubling})$ and $0<D(\Theta)<2^{n/q}$, then for any given $\sigma>0$ and any ball $B(x_0,r)\subset\mathbb R^n$, there exists a constant $C>0$ independent of $f$, $B(x_0,r)$ and $\sigma$ such that
\begin{equation*}
\begin{split}
&\left(\frac{1}{\Theta^q(r)}\cdot\big|\big\{x\in B(x_0,r):\big|[b,M_\alpha](f)(x)\big|>\sigma\big\}\big|\right)^{1/q}\\
\leq& C\cdot\sup_{r>0}\frac{1}{\Theta(r)}\int_{B(x_0,r)}\Phi\left(\frac{|f(x)|}{\sigma}\right)dx,
\end{split}
\end{equation*}
where $\Phi(t)=t(1+\log^+t)$.
\end{corollary}

\end{document}